\documentclass[leqno,oneeqnum,onefignum,onetabnum]{siamltex1213}

\usepackage{amsmath,amssymb}
\usepackage[sc]{mathpazo}
\usepackage{algorithm,algorithmic}
\usepackage{graphicx}
\usepackage{multirow}
\usepackage{makecell}
\usepackage{textcomp}
\usepackage{rotating}
\usepackage{siunitx}
\usepackage{paralist}
\usepackage{booktabs}
\usepackage{lscape}
\usepackage{moresize}
\usepackage{booktabs}
\usepackage[table]{xcolor}

\newcounter{runidnum}
\newcommand{\runid}{\stepcounter{runidnum}\#\therunidnum}
\newcolumntype{g}{>{\columncolor{gray!20}}r}
\usepackage{geometry}
\graphicspath{{figures/}{tables/}}
\makeatletter
\@ifclassloaded{standalone}{}
{
\hypersetup
{
pdfauthor={Andreas Mang and George Biros},
linkcolor=red,
citecolor=orange
}
}
\makeatother%

\newcommand{\figref}[1]{Fig.~\ref{#1}}
\newcommand{\tabref}[1]{Tab.~\ref{#1}}
\newcommand{\secref}[1]{\S\ref{#1}}
\newcommand{\algref}[1]{Alg.~\ref{#1}}
\newcommand{\runref}[1]{run~\##1}
\newcommand{\mcol}[2]{\multicolumn{#1}{r}{#2}}

\newcommand{\tabadjust}{\centering\scriptsize\renewcommand\arraystretch{1.2}}
\newcommand{\algadjust}{\centering\small\renewcommand\arraystretch{1.2}}
\newcommand{\F}[1]{\ensuremath{\mathcal{#1}}}
\newcommand{\D}[1]{\ensuremath{\mathcal{#1}}}
\newcommand{\ns}[1]{\ensuremath{\mathbf{#1}}}
\newcommand{\fs}[1]{\ensuremath{\mathcal{#1}}}
\newcommand{\idiv}{\ensuremath{\nabla\cdot}}
\newcommand{\igrad}{\ensuremath{\nabla}}
\newcommand{\ilap}{\rotatebox[origin=c]{180}{$\nabla$}}
\newcommand{\icurl}{\ensuremath{\nabla \times}}

\newcommand{\id}{\ensuremath{\operatorname{id}}}
\newcommand{\half}[1]{\frac{#1}{2}}
\renewcommand{\d}[1]{\mathop{}\!\mathrm{d}#1}
\newcommand{\dt}{\d{t}}
\newcommand{\dx}{\d{\vect{x}}}

\newcommand{\p} {\partial}

\newcommand{\vect}[1]{\boldsymbol{#1}} 
\newcommand{\mat}[1]{\boldsymbol{#1}}  
\newcommand{\bipa}{\begin{inparaenum}[(\itshape i\upshape)]}
\newcommand{\eipa}{\end{inparaenum}}
\newcommand{\bipasub}{\begin{inparaenum}[(\itshape a\upshape)]}
\newcommand{\eipasub}{\end{inparaenum}}

\newcommand{\ipoint}[1]{\textit{\textbf{#1}}:}
\newcommand{\iname}[1]{\emph{#1}}

\newcommand{\iquote}[1]{``\emph{#1}''}
\newcommand{\ccomment}[1]{{\hfill\color{gray} $\rhd$\,#1}}

\newcommand{\iom}[1]{\int_{\Omega}#1\dx}

\newcommand{\iut}[1]{\int_0^1#1\dt}

\newcommand{\defeq}{\ensuremath{\mathrel{\mathop:}=}}
\newcommand{\eqdef}{\ensuremath{=\mathrel{\mathop:}}}
\newcommand{\T}{\ensuremath{\mathsf{T}}}

\definecolor{sred}{cmyk}{0.01,0.98,0,0.2} 
\reversemarginpar
\newcommand{\mmargin}[1]{{\marginpar{\em\tiny #1}}}\renewcommand{\mmargin}[1]{}

\def\stitle{CONSTRAINED $H^1$-REGULARIZATION SCHEMES FOR DIFFEOMORPHIC REGISTRATION}

\title{Constrained $H^1$-regularization schemes for diffeomorphic image registration\thanks{\ipoint{Publication} SIAM J. Imaging Sci., 9(3), 1154--1194. (41 pages). \url{http://dx.doi.org/10.1137/15M1010919}; \ipoint{Funding} This material is based upon work supported by AFOSR grants FA9550-12-10484 and FA9550-11-10339; by NSF grant CCF-1337393; by the U.S. Department of Energy, Office of Science, Office of Advanced Scientific Computing Research, Applied Mathematics program under award numbers DE-SC0010518 and DE-SC0009286; by NIH grant 10042242; by DARPA grant W911NF-115-2-0121; and by the Technische Universit\"{a}t M\"{u}nchen---Institute for Advanced Study, funded by the German Excellence Initiative (and the European Union Seventh Framework Programme under grant agreement 291763). Any opinions, findings, and conclusions or recommendations expressed herein are those of the authors and do not necessarily reflect the views of AFOSR, DOE, NIH, DARPA, or NSF.}}

\author{Andreas Mang\footnotemark[2]\and George Biros\thanks{The Institute for Computational Engineering and Sciences, The University of Texas at Austin, Austin, Texas, 78712-0027, US ({\tt andreas@ices.utexas.edu}, {\tt gbiros@acm.org}).}}

\begin{document}

\maketitle

\begin{abstract}
We propose regularization schemes for deformable registration and efficient algorithms for their numerical approximation. We treat image registration as a variational optimal control problem. The deformation map is parametrized by its velocity. Tikhonov regularization ensures well-posedness. Our scheme augments standard smoothness regularization operators based on $H^1$- and $H^2$-seminorms with a constraint on the divergence of the velocity field, which resembles variational formulations for Stokes incompressible flows. In our formulation, we invert for a stationary velocity field and a mass source map. This allows us to explicitly control the compressibility of the deformation map and by that the determinant of the deformation gradient. We also introduce a new regularization scheme that allows us to control shear.

We use a globalized, preconditioned, matrix-free, reduced space (Gauss--)Newton--Krylov scheme for numerical optimization. We exploit variable elimination techniques to reduce the number of unknowns of our system; we only iterate on the reduced space of the velocity field. Our current implementation is limited to the two-dimensional case.

The numerical experiments demonstrate that we can control the determinant of the deformation gradient without compromising registration quality. This additional control allows us to avoid oversmoothing of the deformation map. We also demonstrate that we can promote or penalize shear while controlling the determinant of the deformation gradient.
\end{abstract}

\newcommand{\slugmaster}{\slugger{siims}{xxxx}{xx}{x}{x--x}} 



\pagestyle{myheadings}
\thispagestyle{plain}
\markboth
{
ANDREAS MANG AND GEORGE BIROS
}
{
\stitle
}

\section{Introduction}
\label{s:introduction}

Image registration is a key technology in computer vision and imaging sciences. Applications include surveillance, remote sensing, motion tracking, and medical image analysis. Lucid and concise expositions on image registration can be found in~\cite{Modersitzki:2004a, Modersitzki:2009a, Sotiras:2013a}. The problem of image registration can be stated as follows: Given a \iname{reference image} $m_R : \bar{\Omega} \rightarrow \ns{R}$ and a \iname{template image} $m_T : \bar{\Omega} \rightarrow \ns{R}$ with compact support on $\Omega \subset \ns{R}^d$, $d\in\{2,3\}$, we seek a \emph{plausible} map $\vect{y} : \bar{\Omega} \rightarrow \ns{R}^d$ such that the distance between $m_R$ and $m_T\circ \vect{y}$ is as small as possible; $\bar{\Omega}\defeq\Omega\cup\p\Omega$ denotes the closure of $\Omega$ with boundary $\p\Omega$, and the operator $\circ$ is the function composition. If we use an $L^2$-distance to measure the proximity between $m_R$ and $m_T\circ \vect{y}$ we can formulate image registration as a variational optimization problem
\begin{equation}
\label{e:opt-prob-y}
\min_{\vect{y}}\half{1}\|m_R-m_T\circ\vect{y}\|^2_{L^2(\Omega)} + \half{\beta} \|\vect{y}\|^2_{\fs{Y}}.
\end{equation}

\noindent Deformable registration is an ill-posed, nonlinear, and nonconvex optimization problem---regularization is inevitable. The key idea of regularization is to stably compute a solution to a nearby problem. A variety of regularization schemes have been proposed, for example,~\cite{Broit:1981a, Burger:2013a, Christensen:1996a, Chumchob:2011a, Chumchob:2013a, Droske:2003a, Fischer:2002a,Fischer:2003a,FrohnSchauf:2008a,Henn:2005a, Henn:2006a, Mang:2015a, Mang:2012a}. Regularization is typically based on some Tikhonov functional that is added to the objective, which---in the case of~\eqref{e:opt-prob-y}---is a quadratic norm, the contribution of which is controlled by the weight $\beta>0$. The particular choice of the regularization model depends on the application. This is also true for the measure of the proximity between $m_R$ and $m_T\circ\vect{y}$; different choices can be found in~\cite{Modersitzki:2004a,Modersitzki:2009a,Sotiras:2013a}.

A key requirement in many applications, especially in medical imaging, is that the map $\vect{y}$ is a \iname{diffeomorphism}~\cite{Beg:2005a, Burger:2013a, Dupuis:1998a, Trouve:1998a,  Vercauteren:2009a}, i.e., $\vect{y}$ is a bijection, continuously differentiable, and has a continuously differentiable inverse. Formally, we require that $\det(\igrad\vect{y})\not=0$ for every $\vect{x}\in\Omega$, where $\igrad\vect{y}\in\ns{R}^{d\times d}$ is the Jacobian of the deformation map $\vect{y}$. Under the assumption that $\vect{y}$ is orientation preserving we require that $\det(\igrad\vect{y}) > 0$ for every $\vect{x}\in\Omega$. In practice, we would like to control the distance of $\det(\igrad\vect{y})$ from zero.\footnote{Monitoring $\det(\igrad\vect{y})$ does not guarantee that volume elements do not collapse~\cite{Burger:2013a, Haber:2007a, Mang:2015a}. In practice, we have to monitor geometric properties of the deformed grid cells.} Generally speaking, the type and weight of regularization are selected to drive the optimizer to diffeomorphic maps $\vect{y}$ at reasonable computational cost while enabling a good registration between $m_R$ and $m_T$. In the framework of large deformation diffeomorphic image registration we do not directly invert for the deformation map $\vect{y}$ but for its velocity $\vect{v}$. Broadly speaking, we can distinguish between approaches that invert for stationary~\cite{Arsigny:2006a,Ashburner:2007a,Hernandez:2009a,Lorenzi:2013a,Lorenzi:2013b,Vercauteren:2009a} and those that invert for nonstationary velocity fields~\cite{Beg:2005a,Christensen:1996a,Dupuis:1998a,Trouve:1998a}. The proposed formulation uses a stationary velocity, although in principle the extension to a nonstationary velocity is straightforward. In either case, the search space for $\vect{y}$ is typically restricted to the manifold of diffeomorphisms by specifying adequate smoothness requirements for $\vect{v}$~\cite{Ashburner:2011a,Beg:2005a,Dupuis:1998a,Hernandez:2009a,Trouve:1998a}. However, this smoothness control may result in oversmoothing or may lead to $\det(\igrad\vect{y})\approx 0$ or even $\det(\igrad\vect{y})<0$~\cite{Ashburner:2011a}.

Here, we propose constrained regularization schemes for $\vect{v}$ that allow us to control $\det(\igrad\vect{y})$ and the amount of shear in the deformation map $\vect{y}$. We follow up on~\cite{Mang:2015a} where we introduced numerical schemes for our optimal control based large deformation diffeomorphic image registration formulation (for both stationary and nonstationary velocity fields). In particular, we considered two models---one for \iname{compressible} and one for \iname{incompressible} diffeomorphisms $\vect{y}$.\footnote{Related work on incompressible diffeomorphisms can be found in~\cite{Chen:2011a,Chen:2012a,Hinkle:2009a,Mansi:2011a,Ruhnau:2007a,Saddi:2008a}; see \secref{s:related-work}.} For the incompressible case we hypothesized that fixing $\det(\igrad\vect{y})$ to one yields more well-behaved mappings as compared to plain smoothness regularization. We found that enforcing $\det(\igrad\vect{y})=1$ up to numerical accuracy seems to be a too strong constraint for our formulation to be applicable across a wide range of registration problems. We also found that enforcing $\det(\igrad\vect{y})=1$ can lead to excessive shear in the deformation map.

In the present work, we propose new regularization schemes to address these issues. We introduce a mass source $w:\bar{\Omega}\rightarrow\ns{R}$ as an additional unknown to our variational optimization problem. Conceptually, this is equivalent to replacing the incompressibility constraint by a soft constraint (penalty) on the divergence of $\vect{v}$ (see, e.g.,~\cite{Borzi:2002a}). Our formulation avoids ill-conditioning issues in case we set $\idiv\vect{v}$ to a specified value (e.g., zero). We refer to this scheme as \emph{linear Stokes regularization}. Our hypothesis is that the obtained maps $\vect{y}$ are better behaved (smaller variations of the determinant of the deformation gradient) without compromising registration quality as compared to plain smoothness regularization, e.g., used in~\cite{Ashburner:2011a,Avants:2008a,Avants:2011a,Beg:2005a,Dupuis:1998a,Hernandez:2014a,Trouve:1998a,Thirion:1998a,Vercauteren:2009a,Vialard:2012a}. A similar reasoning can be found in connection with hyperelastic regularization models~\cite{Burger:2013a,Droske:2003a}.

Our overarching goal is to design a biophysically constrained framework for large deformation diffeomorphic image registration. Constraints can range from complicated biophysical priors, such as brain tumor models~\cite{Gholami:2016a,Gooya:2012a,Mang:2012b,Hogea:2008a} or cardiac motion models~\cite{Sundar:2009a}, to---like in the present case---simpler models of (nearly) incompressible tissue. The general idea is to favor diffeomorphic deformation maps that have minimal volume changes without compromising data fidelity. An interesting application for incompressible diffeomorphisms is motion estimation in cardiac imaging~\cite{Bistoquet:2008a,Gorce:1996a,Mansi:2011a,Sundar:2009a}. Here, it is expected (at least for healthy individuals) that the volume of the heart muscle does not vary significantly during a cardiac cycle; the deformation map is incompressible. Other applications for (near-)incompressible diffeomorphic registration include time series of abdominal images of a single individual, i.e., images of the liver or the kidneys. Here, we also expect the tissue to mostly behave like an incompressible material. Notice that our new formulation relaxes the incompressibility constraint---it is possible to compute deformation maps that have large local volume changes; we will demonstrate this experimentally.

We, in addition to that, introduce a new regularization scheme that allows us to promote or penalize shear. This formulation also operates in a near-incompressible regime and is motivated from continuum mechanics~\cite{Paterson:1994a, Ranalli:1995a}. In some registration problems the optimizer might drive us to maps that introduce excessive shear. Our new formulation allows us to penalize shear in order to generate maps that are well behaved, guaranteed to be diffeomorphic, and potentially (near-)incompressible. On the contrary, we can also promote shear using the same formulation by simply changing the value of a single parameter. This may be of interest in applications where we expect sharp interfaces (large shear) in the deformation map. We refer to this scheme as \emph{nonlinear Stokes regularization}. We will see that this formulation is---in the limit---equivalent to total variation regularization~\cite{Chen:2012a,Chumchob:2013a,FrohnSchauf:2008a}.

\subsection{Outline of the Method}
\label{s:outline}

We introduce a pseudo time variable $t>0$ and solve for a \iname{stationary velocity field} $\vect{v}\in\fs{V}$, $\fs{V}\subset L^2(\Omega)^d$, and a \iname{mass source} $w\in\fs{W}$, $\fs{W}\subset L^2(\Omega)$, as follows:
\begin{subequations}
\label{e:opt-prob-v-om}
\begin{equation}
\label{e:objective-om}
\min_{\vect{v}, w}\;\;
\half{1}\|m_R - m_1\|^2_{L^2(\Omega)}
+ \half{\beta_v}\|\vect{v}\|_{\fs{V}}^q
+ \half{\beta_w} \|w\|_{\fs{W}}^2
\end{equation}
\noindent subject to
\begin{align}
\p_t m + \igrad m \cdot \vect{v} & = 0
&& {\rm in}\;\Omega\times(0,1],
\label{e:transport-eq-om}
\\
m - m_T &=0
&& {\rm in}\;\Omega\times\{0\},
\label{e:init-transport-eq-om}
\\
\idiv \vect{v} - w &=0
&& {\rm in} \; \Omega,
\label{e:constraint-divv-om}
\end{align}
\end{subequations}

\noindent and periodic boundary conditions on $\p\Omega$. The state variable $m : \bar{\Omega}\times[0,1]\rightarrow \ns{R}$ in~\eqref{e:transport-eq-om} models the transported intensities of the template image $m_T:\bar{\Omega}\rightarrow \ns{R}$ subjected to the stationary velocity field $\vect{v}$. The deformation map $\vect{y}$ is not computed explicitly.\footnote{The Eulerian deformation map $\vect{y}$ and the deformation gradient $\mat{F}_1 \defeq (\igrad\vect{y})^{-1}$ are computed from $\vect{v}$ (see \secref{s:deformation-gradient}).} Instead, the solution of~\eqref{e:transport-eq-om}, i.e.\ $m_1 \defeq m(\,\cdot\,,t=1)$, $m_1:\bar{\Omega}\rightarrow \ns{R}$, represents $m_T\circ\vect{y}$ in~\eqref{e:opt-prob-y}, where $\vect{y}$ represents an Eulerian deformation map. We, like in~\eqref{e:opt-prob-y}, use an $L^2$-distance to measure the proximity between $m_1$ and $m_R$. The objective in~\eqref{e:objective-om} additionally consists of two regularization models that act on the controls $\vect{v}$ and $w$ with weights $\beta_v$ and $\beta_w$, respectively. We provide more details on the choice for the associated norms and the choices for $q>0$ in \secref{s:problem-formulation}.

We augment the regularization on $\vect{v}$ by a constraint on the divergence of the control $\vect{v}$ in~\eqref{e:constraint-divv-om}. Setting $w$ in~\eqref{e:constraint-divv-om} to zero yields a model of incompressible flow~\cite{Chen:2011a,Chen:2012a,Hinkle:2009a,Mang:2015a,Mansi:2011a,Ruhnau:2007a}. This is equivalent to enforcing $\det(\igrad\vect{y})=1$ up to numerical errors (see \cite[p.~70ff.]{Gurtin:1981a}). We relax this incompressibility constraint by introducing an \emph{unknown} mass source $w$, which is determined by solving~\eqref{e:opt-prob-v-om}.

In \secref{s:optimality-conditions} we will see that the optimality system for~\eqref{e:opt-prob-v-om} is a system of space-time nonlinear multicomponent PDEs for the transported intensities $m$, the velocity field $\vect{v}$, the mass source $w$, and the adjoint variables for the transport and divergence condition. Solving this system poses significant challenges. We follow our former work~\cite{Mang:2015a} and solve for the first-order optimality conditions using a globalized, matrix-free, preconditioned inexact (Gauss--)Newton--Krylov method for the Schur complement of the velocity $\vect{v}$. We first derive the optimality conditions and then discretize using a pseudospectral discretization in space with a Fourier basis.

\subsection{Contributions}
\label{s:contributions}

The main goal of this work is to introduce and put to the test new regularization schemes for large deformation diffeomorphic image registration. We use the solver and numerical techniques we have described in \cite{Mang:2015a} to efficiently solve the associated optimization problem. Our main contribution is the formulation and the derivation of the systems. In particular, the contributions are as follows:

\begin{itemize}
\item We extend existing work on continuum mechanical models for incompressible flow~\cite{Borzi:2002a,Chen:2011a,Chen:2012a,Hinkle:2009a,Mang:2015a,Mansi:2011a,Ruhnau:2007a,Saddi:2008a} by introducing a mass source $w$ into the variational optimization problem. This results in a formulation that is more flexible in that we do not fix $\det(\igrad\vect{y})$ to one.
\item We propose a novel $H^1$-regularization scheme that yields a continuum mechanical model with a viscosity that depends on the strain rate tensor (non-Newtonian fluid). This allows us to explicitly control the resistance of the fluid to shear stress and by that promote (shear thinning fluid) or suppress (shear thickening fluid) large shear in the map $\vect{y}$. We will see that this model has a strong resemblance of total variation regularization~\cite{Chen:2012a,FrohnSchauf:2008a}.
\item By using Lagrange multipliers to control the divergence, our formulation avoids ill-conditioning issues in case $\idiv\vect{v}$ is set to a specified value (for example, $w=0$). Our numerical discretization (pseudospectral) allows us to construct fast solvers for the optimality system; we only iterate on the reduced space of the velocity field $\vect{v}$.
\item We provide second order information for numerical optimization.\footnote{Our globalized, matrix-free, preconditioned Newton--Krylov scheme has been described in~\cite{Mang:2015a}. We do not view the numerical scheme as a major contribution but the formulation and the derivation of the associated optimality conditions. A study of our numerical scheme, which includes a comparison to a preconditioned gradient descent algorithm (i.e., our algorithm uses the reduced gradient in the Sobolev space $\fs{V}$), can be found in~\cite{Mang:2015a}.} Although second order methods have widely been used in traditional, variational registration approaches (see, e.g.,~\cite{Modersitzki:2009a}), there has been little work on the use of Newton-type optimization in the framework of large deformation diffeomorphic image registration~\cite{Ashburner:2011a,Hernandez:2014a,Mang:2015a}. Most work in this area is still based on first order numerical optimization strategies~\cite{Avants:2008a,Avants:2011a,Borzi:2002a,Chen:2011a,Chen:2012a,Hart:2009a,Hernandez:2009a,Lee:2010a,Lee:2011a,Vialard:2012a}.
\item We study the effect of incompressibility and smoothness regularization on the overall registration quality as a function of the regularization parameters. We quantify registration accuracy in terms of overlap measures and compare our results against the diffeomorphic {\tt DEMONS} algorithm~\cite{Vercauteren:2009a}. We demonstrate that, by introducing a constraint on $\idiv\vect{v}$, we can control $\det(\igrad\vect{y})$ without compromising registration quality. We show that our model allows us to avoid oversmoothing of the deformation map. We also study the effect of controlling shear.
\end{itemize}

\subsection{Limitations and Unresolved Issues}
\label{s:limitations}

Here, we summarize the limitations and unresolved issues of our work:
\begin{itemize}
\item We introduce an additional regularization parameter. This makes it more difficult to design a black-box solver and more expensive to automatically calibrate the algorithm.
\item We assume similar intensity statistics of $m_R$ and $m_T$. This is a common assumption in many image registration approaches~\cite{Beg:2005a,Chen:2011a,FrohnSchauf:2008a,Hart:2009a,Hernandez:2009a,Lee:2010a,Museyko:2009a,Vialard:2012a}. For multimodal registration problems different distance measures have to be considered (see e.g.,~\cite{Modersitzki:2004a,Sotiras:2013a}).
\item We present results only in two dimensions. Nothing in our formulation and numerical approximation is specific to the two-dimensional case. Overall, the method is very expensive and a practical three-dimensional implementation requires more work.
\item We only report results for stationary velocities (see, e.g.,~\cite{Ashburner:2007a,Hernandez:2009a,Vercauteren:2009a}). We have implemented and tested time-varying velocities (in~\cite{Mang:2015a} we report results for incompressible $\vect{y}$). For a two-image registration problem we found that a velocity that changes in time does not improve the quality of the registration. For tracking problems like optical flow~\cite{Borzi:2002a,Chen:2011a,Chen:2012a,Horn:1981a,Kalmoun:2011a,Ruhnau:2007a} or time series of medical images~\cite{Mang:2008a} a time-dependent velocity may be necessary. Nothing changes in our formulation, just the problem size (see~\cite{Mang:2015a}).
\end{itemize}

\subsection{Related Work}
\label{s:related-work}

There is a vast body of literature on image registration. Here, we restrict the discussion to approaches that are closely related to our work. We refer to~\cite{Modersitzki:2004a,Modersitzki:2009a,Sotiras:2013a} for a more general overview on algorithmic developments and formulations.

Our approach shares numerous characteristics with methods that have appeared in the past. Optimal control formulations for image registration have been discussed in~\cite{Borzi:2002a,Chen:2011a,Hart:2009a,Lee:2010a,Lee:2011a,Mang:2015a,Vialard:2012a}. Our work is related to large deformation diffeomorphic metric mapping~\cite{Avants:2008a,Avants:2011a,Beg:2005a,Dupuis:1998a,Trouve:1998a,Younes:2007a} (see~\cite{Mang:2015a} for a connection). It differs in that we, like~\cite{Arsigny:2006a,Ashburner:2007a,Hernandez:2009a,Lorenzi:2013a,Lorenzi:2013b,Vercauteren:2009a}, invert for a stationary velocity field. Our formulation also shares conceptual ideas with traditional optical flow formulations~\cite{Horn:1981a,Kalmoun:2011a,Ruhnau:2007a}. We refer to~\cite{Mang:2015a} for a more detailed discussion. In this review we will focus on approaches that \bipa\item introduce mass conservation as an additional constraint and/or \item aim at recovering motion fields that locally contain significant shear\eipa.

One way to explicitly control $\det(\igrad\vect{y})$ is to set it to one. This is equivalent to working with incompressible velocity fields (see \cite[p.~77ff.]{Gurtin:1981a}); we refer to this model as \iname{linear Stokes regularization}. Formulations based on divergence-free velocity fields have been described in~\cite{Borzi:2002a,Chen:2011a,Hinkle:2009a,Mang:2015a,Mansi:2011a,Ruhnau:2007a,Saddi:2008a}. None of these consider an inversion for a mass source $w$. All of these, with the exception of our preceding work~\cite{Mang:2015a}, use first order information only for numerical optimization. Other formulations for controlling $\det(\igrad\vect{y})$ can be found in~\cite{Aganj:2015a,Baluwala:2013a,Burger:2013a,Haber:2010a,Haber:2004a,Haber:2007a,Li:2009a,Loeckx:2004a,Modersitzki:2008a,Nielsen:2002a,Pennec:2005a,Rohlfing:2003a,Ruehaak:2013a,Sdika:2008a,Yanovsky:2007a}.

We, in addition to that, introduce a continuum mechanical model that controls shear (either promoting or penalizing it). We refer to this formulation as \iname{nonlinear Stokes regularization}. Related approaches based on nonquadratic regularization norms ($L^1$-norm or total variation) have been described in~\cite{Chen:2012a,Chumchob:2013a,FrohnSchauf:2008a,Ruan:2009a,Zach:2007a}. In our formulation, the regularization is a function of the shear strain rate. This allows us to explicitly control the amount of shear in the deformation map. Whether such an explicit control is beneficial in certain applications remains to be seen. We will see that our formulation is in the limit equivalent to total variation regularization~\cite{Chen:2012a,Chumchob:2013a,FrohnSchauf:2008a}. Our model couples the individual components of the regularized vector field as opposed to component wise vectorial total variation regularization~\cite{Chen:2012a,FrohnSchauf:2008a}.

Other approaches for estimating an expected discontinuous motion field include locally adaptive (i.e., direction-dependent and/or intensity-driven) regularization~\cite{Pace:2013a,Papiez:2014a,SchmidtRichberg:2012a}, are based on a decomposition of the bodyforce~\cite{Baluwala:2013a}, or are based on a subdivision of the domain~\cite{Risser:2013a,Ruehaak:2013a,Wu:2008a}. The formulations in~\cite{Chumchob:2013a,Chumchob:2011a,FrohnSchauf:2008a,Ruan:2009a,SchmidtRichberg:2012a,Pace:2013a} operate on the deformation map or the displacement field. Our formulation operates on the velocity field instead and as such falls into the category of large deformation models. We can, like in the linear Stokes case, control the magnitude of $\det(\igrad\vect{y})$. All mentioned approaches for estimating sliding motion, with the exception of~\cite{Baluwala:2013a,Ruan:2009a,Ruehaak:2013a}, do not feature such a control. Our formulation---as opposed to~\cite{Baluwala:2013a,Pace:2013a,Risser:2013a,Ruehaak:2013a,SchmidtRichberg:2012a,Wu:2008a}---does not require any partitioning of the domain (presegmentation). We exploit the second order Newton--Krylov scheme we have introduced in~\cite{Mang:2015a}, which further distinguishes us from most of the preceding work. Our formulation allows us to control the shear within the estimated motion field on the basis of a single, strictly positive parameter; we can not only promote shear but also penalize it.

\subsection{Organization and Notation}
\label{s:organization}

We provide additional details on the optimal control formulation in~\secref{s:problem-formulation}. The optimality system is summarized in \secref{s:optimality-conditions}. The numerics are described in \secref{s:numerics} (we refer to~\cite{Mang:2015a}, where we originally introduced our numerical scheme, for more details). We report experiments in \secref{s:numerical-experiments} and conclude with \secref{s:conclusions}. Additional derivations, comments, algorithmic details, and measures of registration performance can be found in the appendix.

An overview of the commonly used symbols can be found in \tabref{t:notation}. Vectorial quantities and matrices are denoted in boldface. Function spaces, differential operators, and functionals are denoted in calligraphy. A superscript $h$ is added to the variables whenever we refer to discretized quantities.

\begin{table}[t]
\centering
\caption{Commonly used notation and symbols.}
\label{t:notation}
\begin{small}
\begin{tabular}{ll}
\toprule
Symbol/Notation
&
Description
\\\midrule
CFL         & Courant--Friedrichs--Lewy (condition)\\
DSC         & Dice similarity coefficient\\
FFT         & fast Fourier transform\\
FNE         & false negative error\\
FPE         & false positive error\\
JSC         & Jaccard similarity coefficient \\
KKT         & Karush--Kuhn--Tucker\\
matvec      & matrix-vector product\\
PDE         & partial differential equation\\
PDE solve   & solution of the hyperbolic transport equations\\
PCG         & preconditioned conjugate gradient (method)\\
RK2         & 2nd order Runge--Kutta (method)\\
\midrule
$d$         & spatial dimensionality; typically $d\in\{2,3\}$\\
$\Omega$    & spatial domain; $\Omega\defeq(-\pi,\pi)^d\subset\ns{R}^d$ with boundary $\p\Omega$ and closure $\bar{\Omega}\defeq\Omega\cup\p\Omega$\\
$\vect{x}$  & spatial coordinate; $\vect{x}\defeq(x^1,\ldots,x^d)^\T\in\ns{R}^d$ \\
$m_R$       & reference image; $m_R : \bar{\Omega} \rightarrow \ns{R}$\\
$m_T$       & template image; $m_T : \bar{\Omega} \rightarrow \ns{R}$\\
$m$         & state variable (transported intensities); $m : \bar{\Omega}\times[0,1]\rightarrow\ns{R}$\\
$m_1$       & deformed template image (state variable at $t=1$); $m_1 : \bar{\Omega} \rightarrow \ns{R}$\\
$\lambda$   & adjoint variable (transport equation); $\lambda : \bar{\Omega}\times[0,1]\rightarrow\ns{R}$\\
$p$         & adjoint variable (incompressibility constraint); $p : \bar{\Omega}\rightarrow\ns{R}$\\
$\vect{v}$  & control variable (stationary velocity field); $\vect{v} : \bar{\Omega} \rightarrow \ns{R}^d$\\
$w$         & control variable (mass source); $w : \bar{\Omega} \rightarrow \ns{R}$\\
$\vect{b}$  & bodyforce; $\vect{b} : \bar{\Omega} \rightarrow \ns{R}^d$\\
$\D{H}$     & (reduced) Hessian\\
$\vect{g}$  & (reduced) gradient\\
$\mat{F}_1$ & deformation gradient at $t=1$; $\mat{F}_1 : \bar{\Omega}\rightarrow\ns{R}^{d\times d}$; $\mat{F}_1 = \mat{F}(\cdot,t=1)$, $\mat{F}_1\equiv(\igrad\vect{y})^{-1}$\\
$\beta_v$   & regularization parameter for the control $\vect{v}$\\
$\beta_w$   & regularization parameter for the control $w$\\
$\nu$       & Glen's flow law exponent; $\nu > 0$\\
$\D{E}$     & strain rate tensor; $\D{E}[\vect{v}]\defeq\half{1}((\igrad\vect{v})+(\igrad\vect{v})^\T)$\\
$\D{A}$     & regularization operator (variation of regularization model acting on $\vect{v}$)\\
$\igrad$    & gradient operator (acts on scalar and vector fields)\\
$\ilap$     & Laplacian operator (acts on scalar and vector fields)\\
$\idiv$     & divergence operator (acts on vector and 2nd order tensor fields)
\\\bottomrule
\end{tabular}
\end{small}
\end{table}

\section{Problem Formulation}
\label{s:problem-formulation}

The images to be registered are modeled as compactly supported functions on the domain $\Omega\defeq(-\pi,\pi)^d\subset\ns{R}^d$, $d\in\{2,3\}$, with boundary $\p\Omega$, and closure $\bar{\Omega}\defeq \Omega \cup \p\Omega$. We introduce a pseudo time variable $t>0$ and solve for a stationary velocity field $\vect{v}\in\fs{V}$, $\fs{V}\subset L^2(\Omega)^d$, and a mass source $w\in\fs{W}$, $\fs{W} \subset L^2(\Omega)$, as follows:
\begin{subequations}
\label{e:opt-prob-v}
\begin{equation}
\label{e:objective}
\min_{\vect{v},w}
\F{J}[\vect{v},w]
\defeq
\half{1}\|m_R - m_1\|_{L^2(\Omega)}^2
+ \half{\beta_v} \|\vect{v}\|^q_{\fs{V}}
+ \gamma\half{\beta_w} \|w\|^2_{\fs{W}}
\end{equation}
\noindent subject to
\begin{align}
\p_t m + \igrad m \cdot \vect{v}
& = 0
&& {\rm in}\;\Omega\times(0,1],
\label{e:transport-eq}
\\
m - m_T
&=0
&& {\rm in}\;\Omega\times\{0\},
\label{e:init-transport-eq}
\\
\gamma(\idiv \vect{v} - w)
&=0
&& {\rm in} \; \Omega,
\label{e:constraint-divv}
\end{align}
\end{subequations}

\noindent and periodic boundary conditions on $\p\Omega$. We do not directly invert for the map $\vect{y}$ but for its velocity $\vect{v}$. This is different from the problem formulation in \eqref{e:opt-prob-y}; in our formulation, the solution of~\eqref{e:transport-eq}---$m_1 \defeq m(\cdot, t=1)$, $m_1:\bar{\Omega}\rightarrow\ns{R}$---is equivalent to $m_T\circ\vect{y}$ in~\eqref{e:opt-prob-y}, where $\vect{y}$ is the deformation map defined in an Eulerian frame of reference. We, as in~\eqref{e:opt-prob-y}, use an $L^2$-distance to measure the proximity between the deformed template image $m_1$ and the reference image $m_R$. This is a common choice in many deformable image registration algorithms (see, e.g.,~\cite{Beg:2005a, Lee:2010a, Museyko:2009a, Vialard:2012a}). The parameters $\beta_v>0$ and $\beta_w>0$ control the contribution of the regularization norms. The parameter $\gamma\in\{0,1\}$ is introduced for clarity. If we set $\gamma=0$ we obtain a formulation that is related to available models for large deformation diffeomorphic image registration~\cite{Ashburner:2007a,Beg:2005a,Hart:2009a,Hernandez:2009a,Vialard:2012a,Vercauteren:2009a} (see~\cite{Hart:2009a, Mang:2015a} for a more detailed insight).\footnote{We invert for a stationary velocity field. Other work on diffeomorphic registration based on stationary velocity fields can, for instance, be found in~\cite{Arsigny:2006a,Ashburner:2007a,Hernandez:2009a,Lorenzi:2013a,Lorenzi:2013b,Vercauteren:2009a}. The traditional large deformation diffeomorphic metric mapping framework introduced in~\cite{Beg:2005a,Dupuis:1998a,Trouve:1995a,Trouve:1998a,Younes:2007a} uses time-dependent velocity fields. Our implementation allows us to invert for time-dependent and stationary velocity fields alike (in case quadratic regularization models are considered); nothing changes in our formulation, just the number of unknowns (see~\cite{Mang:2015a} for details). Here, we limit ourselves to stationary $\vect{v}$.} If we set $\gamma$ to one and $w$ to zero, we obtain a model for incompressible diffeomorphisms; i.e., we enforce $\det(\mat{F}_1) = 1$ up to numerical accuracy; the tensor field $\mat{F}_1:\bar{\Omega}\rightarrow\ns{R}^{d\times d}$ is the deformation gradient at $t=1$ computed from $\vect{v}$ (see \secref{s:deformation-gradient} for details). Similar approaches for incompressible diffeomorphisms have been described in~\cite{Chen:2011a, Hinkle:2009a, Mang:2015a, Mansi:2011a, Ruhnau:2007a, Saddi:2008a}. We extend these by introducing a nonzero mass source $w$. This allows us to relax the model from incompressible diffeomorphisms to a model of near-incompressible diffeomorphisms. The regularization on $w$ in~\eqref{e:objective} acts like a penalty on $\idiv\vect{v}$; we use an $H^1$-norm:
\begin{equation}
\label{e:reg-w}
\|w\|_{\fs{W}}^2 = \|w\|^2_{H^1(\Omega)}
=\iom{\igrad w \cdot \igrad w + w^2}.
\end{equation}

\noindent We use $H^1$- and $H^2$-seminorms to regularize $\vect{v}$; in particular,
\begin{equation}
\label{e:reg-v}
|\vect{v}|^2_{H^1(\Omega)}
\defeq\iom{\igrad\vect{v}:\igrad\vect{v}}
\quad\text{and}\quad
|\vect{v}|^2_{H^2(\Omega)}
\defeq\iom{\ilap\vect{v}\cdot\ilap\vect{v}},
\end{equation}

\noindent respectively. The choice of an $H^2$-seminorm is motivated by related work on large deformation diffeomorphic image registration (see, e.g.,~\cite{Beg:2005a,Hart:2009a,Vialard:2012a}). The choice of an $H^1$-seminorm is motivated by the fact that we obtain optimality conditions that are similar to Stokes equations in fluid mechanics; we refer to this (near-)incompressible formulation as \iquote{linear Stokes regularization}.

Since we observed that a model of incompressible flow may promote shear, we additionally introduce a nonlinear regularization model that allows us to control (promote or penalize) shear in the deformation field in a problem-dependent way. This model is motivated from continuum mechanics\footnote{We will arrive at a formulation that resembles continuum mechanical models that can be found in geoscience applications~\cite{Isaac:2015a, Paterson:1994a, Ranalli:1995a}.} and given by
\begin{equation}
\label{e:shear-control-reg-v}
|\vect{v}|_{H^1(\Omega)}^{(1+\nu)/2 \nu}
=
\frac{2\nu}{\nu+1}
\iom
{
\left(
\half{1} ((\igrad \vect{v}) + (\igrad\vect{v})^\T)
:
\half{1} ((\igrad \vect{v}) + (\igrad\vect{v})^\T)
\right)^{(1+\nu)/2 \nu}
}.
\end{equation}

\noindent Here, $\nu>0$ controls the nonlinearity and $\half{1} ((\igrad \vect{v}) + (\igrad\vect{v})^\T)\eqdef\D{E}[\vect{v}]$ is the \iname{strain rate tensor}. We will see that we arrive at a Stokes-like optimality system with a viscosity that depends on the strain rate (see \secref{s:optimality-conditions} for details). For $\nu \in(0,1)$ we obtain a \iname{shear thickening} and for $\nu > 1$ a \iname{shearthinning} fluid. Notice that we approach a total variation like regularization model as $\nu$ in~\eqref{e:shear-control-reg-v} tends to $\infty$ (see \secref{s:total-variation}). Thus, we can explicitly control the shear within $\vect{y}$ via $\nu$. This model, in combination with the incompressibility constraint, yields a deformation map for which $\det(\mat{F}_1)=1$. This is a fundamental difference to most existing models for estimating sliding motion\footnote{We note that our formulation allows us to only approximate discontinuous motion fields (sliding motion) by promoting shear; the computed map will still be continuous.} (with the exception of~\cite{Baluwala:2013a, Ruan:2009a, Ruehaak:2013a}), since these in general do not explicitly control the determinant of the deformation gradient. We have also tested a version of this model with a relaxed incompressibility constraint. We refer to this formulation as \iquote{nonlinear Stokes regularization}. Notice, that the derivation we describe in the present work will also hold if we replace \eqref{e:shear-control-reg-v} with a total variation regularization model.

\textbf{Remark on some theoretical considerations}: There are several questions at hand. A first question regards an appropriate choice of the space for $\vect{v}$ so that the transport equation~\eqref{e:transport-eq} has a unique solution and preserves the smoothness of the initial image $m(\cdot,0) = m_T$. The answer to this question depends on smoothness of the input images $m_T$ and $m_R$. A second question regards the sufficient regularity of the adjoint variables, which is required to justify a gradient-descent scheme for solving~\eqref{e:opt-prob-v}. Finally, a third question regards the existence and uniqueness of the solution of~\eqref{e:opt-prob-v}. If the input images $m_T$ and $m_R$ are adequately smooth, and the regularization space and weights are large enough, then the transport equation has a unique solution, smoothness is preserved, the adjoint variables are smooth, and the problem has a solution. Uniqueness requires an even stronger regularization to ensure the convexity of the problem. In our experiments we always use smooth images (we apply a Gaussian smoothing operator to the input images). However, it is not known if our regularization scheme for the velocity is in the theoretical limit sufficient for all three questions to have an affirmative answer. Numerically, we control the velocity by adjusting the regularization weight to ensure we obtain diffeomorphic maps. Informally, our experimental studies suggest that $H^1$-regularity of $\vect{v}$ and the penalization of $\idiv\vect{v}$ can provide sufficient smoothness as long as the regularization parameters are sufficiently large. We provide a lengthier discussion in \secref{s:theoretical-considerations}, based on work of other authors~\cite{Chen:2011b,Chen:2011a,Crippa:2007a,DiPerna:1989a}. A detailed theoretical analysis is beyond the scope of this work and remains open for future work.

\section{Optimality Conditions}
\label{s:optimality-conditions}

We use the method of Lagrange multipliers to solve~\eqref{e:opt-prob-v}. The Lagrangian reads
\begin{align}
\label{e:lagrangian}
\F{L}[\vect{\phi}]
\defeq
&
\F{J}[\vect{v},w]
+
\int_0^1
\langle
\p_t m + \igrad m \cdot \vect{v},\lambda
\rangle_{L^2(\Omega)}
\d{t}
+
\langle m - m_T,\mu \rangle_{L^2(\Omega)}
-
\langle \idiv \vect{v} - w,p\rangle_{L^2(\Omega)}
\end{align}

\noindent with Lagrange multipliers $\lambda:\bar{\Omega}\times[0,1]\rightarrow\ns{R}$ for the hyperbolic transport equation~\eqref{e:transport-eq}, $\mu:\bar{\Omega}\rightarrow\ns{R}$ for the initial condition~\eqref{e:init-transport-eq}, and $p:\bar{\Omega}\rightarrow\ns{R}$ for the incompressibility constraint~\eqref{e:constraint-divv} ($p$ is referred to as \iname{pressure} in fluid dynamics); $\vect{\phi}\defeq(m,\vect{v},w,\lambda,\mu,p)$ and $\langle \cdot,\cdot\rangle_{L^2(\Omega)}$ denotes the standard $L^2$ inner product defined on $\Omega$.

Our algorithm falls into the class of reduced space Newton--Krylov methods~\cite{Biros:2005a,Biros:2005b}. This will also be reflected by the optimality systems we present below. The interested reader is referred to \secref{s:algorithm} for a more detailed explanation of the optimality systems, the conceptual ideas behind our algorithm, and details for its implementation. We describe our Newton--Krylov algorithm in more detail in~\cite{Mang:2015a}. We refer to~\cite{Borzi:2012a,Gunzburger:2003a} for general information on optimal control theory and PDE constrained optimization; the conceptual ideas we use for our optimization scheme are described in~\cite{Nocedal:2006a}. We use an \iname{optimize-then-discretize} approach.\footnote{For a discussion on advantages and disadvantages of the optimize-then-discretize and the discretize-then-optimize approach we refer to~\cite{Gunzburger:2003a}.} The resulting optimality conditions is what we discuss next.

\subsection{First order Optimality System}
\label{s:optimality-system}

From Lagrange multiplier theory we know that the variations of $\F{L}$ with respect to all variables have to vanish for an admissible solution of~\eqref{e:opt-prob-v}. Taking variations of $\F{L}$ in~\eqref{e:lagrangian} with respect to $m$, $\vect{v}$, $w$, $\lambda$, $\mu$, and $p$, in directions $\tilde{m}$, $\vect{\tilde{v}}$, $\tilde{w}$, $\tilde{\lambda}$, $\tilde{\mu}$, and $\tilde{p}$, and applying integration by parts, yields the optimality system (i.e., the \emph{first order necessary optimality conditions} (KKT conditions) in strong form)
\begin{subequations}
\label{e:first-order-opt}
\begin{align}
\p_t m + \igrad m \cdot \vect{v} & = 0
&&{\rm in}\;\; \Omega \times (0,1],
\label{e:state-eq}
\\
m  - m_T &= 0
&&{\rm in}\;\; \Omega \times\{0\},
\label{e:init-state-eq}
\\
-\p_t \lambda - \idiv (\vect{v}\lambda) & = 0
&&\text{in}\;\; \Omega \times [0,1),
\label{e:adj-eq}
\\
\lambda + (m  - m_R) &= 0
&&{\rm in}\;\; \Omega \times\{1\},
\label{e:terminal-adj-eq}
\\
\gamma(\idiv\vect{v} - w)  & = 0
&&\text{in}\;\;\Omega,
\label{e:ic-constraint}
\\
\vect{g}_v\defeq
\beta_v\D{A}[\vect{v}] + \gamma\igrad p + \vect{b} & = 0
&&\text{in}\;\; \Omega,
\label{e:control-eq-v}
\\
g_w\defeq
\gamma(\beta_w(-\ilap w + w) + p) & = 0
&&\text{in}\;\; \Omega,
\label{e:control-eq-w}
\end{align}
\end{subequations}

\noindent subject to periodic boundary conditions on $\p\Omega$. The parameter $\gamma\in\{0,1\}$ enables or disables the constraint on the divergence of $\vect{v}$ in \eqref{e:ic-constraint}. Further, \[\vect{b}\defeq \iut{\lambda\igrad m}\] is the \iname{body force}. The operator $-\ilap + \id$ (where $\id$ is the identity operator) in~\eqref{e:control-eq-w} is the first variation of the $H^1$-norm in~\eqref{e:reg-w}. The operator $\D{A}$ in~\eqref{e:control-eq-v} is the first variation of the regularization model for $\vect{v}$. In particular, we have
\begin{equation}
\label{e:var-reg-v}
\D{A}[\vect{v}] = -\ilap \vect{v}
\qquad
\text{and}
\qquad
\D{A}[\vect{v}] = \ilap^2\vect{v}
\end{equation}

\noindent for the $H^1$- and the $H^2$-seminorm in~\eqref{e:reg-v}, respectively.\footnote{For $\D{A}=-\ilap$, $\gamma = 1$, and $w=0$ we obtain a linear Stokes regularization model (i.e., a model for incompressible diffeomorphisms; see~\cite{Mang:2015a}).} Further, we have
\begin{equation}
\label{e:1st-var-reg-v-nls}
\D{A}[\vect{v}] = - \idiv2\eta[\vect{v}] \D{E}[\vect{v}]
\end{equation}

\noindent if we consider the regularization model in~\eqref{e:shear-control-reg-v}, where $\eta[\vect{v}] \defeq\operatorname{tr}(\D{E}[\vect{v}]\D{E}[\vect{v}])^{(1-\nu)/2\nu}$ is the effective viscosity, $\nu>0$ is Glen's flow law exponent, and $\D{E}[\vect{v}] \defeq \half{1}((\igrad\vect{v})+(\igrad\vect{v})^\T)$ is the strain rate tensor.

In the language of optimal control~\eqref{e:state-eq} is referred to as the \iname{state equation} (with initial condition~\eqref{e:init-state-eq}), \eqref{e:adj-eq} as the \iname{adjoint equation} (with final condition~\eqref{e:terminal-adj-eq}), and~\eqref{e:control-eq-v} and~\eqref{e:control-eq-w} as the \iname{control equations}, respectively. Notice that the adjoint equation models the transport of the mismatch between $m_1$ (deformed template image) and $m_R$ backward in time; $\lambda$ will (ideally) tend to zero if we approach the solution of our problem.

We can directly use the optimality system in~\eqref{e:first-order-opt} to design an iterative scheme for computing a solution to~\eqref{e:opt-prob-v}. This will result in a first order gradient descent scheme, which is still widely used in large deformation diffeomorphic image registration~\cite{Beg:2005a,Hart:2009a,Vialard:2012a} despite its linear convergence. As we have seen in~\cite{Mang:2015a} (for the compressible and the incompressible case) preconditioned gradient descent schemes\footnote{The control equation provides the reduced $L^2$ gradient (variation of the objective with respect to $\vect{v}$); by \emph{preconditioned} gradient descent we mean that we use the gradient in the Sobolev space $\fs{V}$ for numerical optimization.} are inferior to preconditioned Newton--Krylov schemes in case we strive for high inversion accuracy. However, exploiting second order information requires more work; we have to derive the second variations of the Lagrangian. This is what we present next.

\subsection{Newton Step}
\label{s:newton-step}

We use a globalized, inexact, reduced space \mbox{(Gauss--)}Newton--Krylov method for numerical optimization (see~\secref{s:numerics}); we solve~\eqref{e:first-order-opt} using a Newton linearization. We have to compute variations of the weak form of the optimality conditions in~\eqref{e:first-order-opt}, i.e., second variations of the Lagrangian in \eqref{e:lagrangian}, to obtain the associated KKT system. Following the standard theory of calculus of variations, invoking the appropriate Green's identities (integration by parts), and adhering to the fact that we consider a reduced space method, we arrive at the following system (which corresponds to the strong form of the second variations of the Lagrangian in \eqref{e:lagrangian}):
\begin{subequations}
\label{e:newton-step}
\begin{align}
\p_t \tilde{m} + \igrad \tilde{m} \cdot \vect{v}
+ \igrad m \cdot \tilde{\vect{v}} & = 0
&&\text{in}\;\; \Omega \times (0,1],
\label{e:inc-state-eq}
\\
\tilde{m} &= 0
&&\text{in}\;\;\Omega \times\{0\},
\label{e:init-inc-state-eq}
\\
-\p_t \tilde{\lambda} - \idiv(\tilde{\lambda}\vect{v})
-\idiv(\lambda\tilde{\vect{v}}) & = 0
&&\text{in}\;\; \Omega \times [0,1),
\label{e:inc-adj-eq}
\\
\tilde{\lambda} +\tilde{m} &= 0
&&\text{in}\;\;\Omega \times\{1\},
\label{e:terminal-inc-adj-eq}
\\
\gamma(\idiv \tilde{\vect{v}} - \tilde{w}) & = 0
&&\text{in}\;\; \Omega,
\label{e:inc-ic-constraint}
\\
\beta_v\D{B}[\tilde{\vect{v}}]
+ \gamma\igrad\tilde{p} + \vect{\tilde{b}}
& = -\vect{g}_v
&&\text{in}\;\;\Omega,
\label{e:inc-control-eq-v}
\\
\gamma(\beta_w(-\ilap\tilde{w} + \tilde{w}) + \tilde{p}) & = -g_w
&&\text{in}\;\; \Omega,
\label{e:inc-control-eq-w}
\end{align}
\end{subequations}

\noindent with periodic boundary conditions on $\p\Omega$ and \iname{incremental body force}
\[
\vect{\tilde{b}} =
\iut{\tilde{\lambda} \igrad m
+ \lambda \igrad \tilde{m}}.
\]

\noindent We refer to~\eqref{e:inc-state-eq} (with initial condition~\eqref{e:init-inc-state-eq}), \eqref{e:inc-adj-eq} (with final condition~\eqref{e:terminal-inc-adj-eq}), \eqref{e:inc-control-eq-v} and~\eqref{e:inc-control-eq-w} as \iname{incremental state}, \iname{incremental adjoint}, and \iname{incremental control equations}, respectively. The incremental variables are denoted with a tilde. The incremental control equations~\eqref{e:inc-control-eq-w} and~\eqref{e:inc-control-eq-v} represent the action of the reduced space Hessian operators on the control variables (Hessian matvec). The right-hand sides in~\eqref{e:inc-ic-constraint} and~\eqref{e:inc-control-eq-w} correspond to the reduced gradients in~\eqref{e:control-eq-w} and~\eqref{e:control-eq-v}, respectively.

The operator $\D{B}$ is the second variation of the regularization model acting on $\vect{v}$. It coincides with the first variations in~\eqref{e:var-reg-v} for the quadratic regularization models in~\eqref{e:reg-v}. This also holds for the second variation of the $H^1$-norm in~\eqref{e:reg-w} (see~\eqref{e:inc-control-eq-w}). The second variation for the nonlinear regularization model in~\eqref{e:shear-control-reg-v} does not coincide with its first variation; we obtain
\begin{equation}
\label{e:2nd-var-reg-v-nls}
\D{B}(\vect{v})[\vect{\tilde{v}}]
=
-\idiv 2\eta[\vect{v}]
\bigg(
\D{I}
+
\underbrace
{
\frac{1-\upsilon}{\upsilon}
\frac{\D{E}[\vect{v}] \otimes \D{E}[\vect{v}]}
{\D{E}[\vect{v}]:\D{E}[\vect{v}]}
}_{\eqdef\D{Q}[\vect{v}]}
\bigg)
\D{E}[\tilde{\vect{v}}]
\end{equation}

\noindent instead, where $\otimes$ is the tensor outer product and $\F{I}$ is a fourth order identity tensor.

We can significantly simplify these systems by exploiting variable elimination techniques. This allows us to merely iterate on the reduced space of the velocity field $\vect{v}$. This is what we discuss next.

\subsection{Reduced Systems}
\label{s:elimination}

We eliminate the control and adjoint variables $w$ and $p$, and by that the constraint on the divergence of the velocity field $\vect{v}$ from the optimality system~\eqref{e:first-order-opt}. The systems we provide below are the ones we solve numerically (see \secref{s:algorithm}). We provide details on their derivation in \secref{s:constraint-elimination-derivation}. We arrive at
\begin{subequations}
\label{e:first-order-opt-elim}
\begin{align}
\p_t m + \igrad m \cdot \vect{v} & = 0
&&{\rm in} \;\; \Omega \times (0,1],
\label{e:state-eq-elim}
\\
m - m_T &= 0
&&{\rm in} \;\; \Omega \times\{0\},
\\
-\p_t\lambda - \idiv (\lambda\vect{v}) & = 0
&&\text{in}\;\;\Omega\times[0,1),
\label{e:adj-eq-elim}
\\
\lambda + (m - m_R) &= 0
&&{\rm in} \;\; \Omega \times\{1\},
\label{e:init-adj-eq-elim}
\\
\vect{g} \defeq \beta_v\D{A}[\vect{v}] + \D{K}[\vect{b}] & = 0
&&\text{in}\;\;\Omega,
\label{e:control-eq-elim}
\end{align}
\end{subequations}

\noindent with periodic boundary conditions on $\p\Omega$ to replace~\eqref{e:first-order-opt}. The operator $\D{A}$ corresponds to the first variation of the regularization models. The operator $\D{K}$ projects $\vect{v}$ onto the space of near-incompressible velocity fields. If we consider the regularization models in~\eqref{e:reg-v} we have
\begin{equation}
\label{e:project-body-force}
\D{K}[\vect{b}]
=
-
\igrad
\D{M}^{-1}
\ilap^{-1}
\idiv\vect{b}
+ \vect{b},
\quad
\text{where}
\;
\D{M}=\beta_v(\beta_w(-\ilap+\id))^{-1}+\id,
\end{equation}
\noindent and $\id$ is the identity operator. If we set $w=0$ this operator simplifies to $\D{M}=\id$. If we consider \eqref{e:shear-control-reg-v} instead, we have
\[
\D{K}[\vect{b},\vect{v}]
=
\igrad
\D{M}^{-1}
\ilap^{-1}
\idiv
(
\idiv2\beta_v\hat{\eta}[\vect{v}]\D{E}[\vect{v}]-\vect{b}
)
+\vect{b},
\]

\noindent where $\D{M}=2\beta_v\bar{\eta}[\vect{v}](\beta_w(-\ilap+\id))^{-1}+\id$. If we set $w=0$ we obtain $\D{M} = \id$.

The system no longer depends on $w$ and $p$. This allows us to efficiently solve~\eqref{e:opt-prob-v}; we only iterate on the reduced space of the velocity field $\vect{v}$. Computing variations of the weak form of~\eqref{e:first-order-opt-elim} yields the Newton step
\begin{subequations}
\label{e:second-order-elim}
\begin{align}
\p_t \tilde{m} + \igrad \tilde{m} \cdot \vect{v}
+ \igrad m \cdot \vect{\tilde{v}} & = 0
&&\text{in}\;\; \Omega \times (0,1],
\label{e:inc-state-eq-elim}
\\
\tilde{m} &= 0
&&\text{in}\;\;\Omega\times\{0\},
\\
-\p_t \tilde{\lambda} - \idiv(\tilde{\lambda}\vect{v})
-\idiv(\lambda\vect{\tilde{v}}) & = 0
&&\text{in}\;\; \Omega \times [0,1),
\label{e:inc-adj-eq-elim}
\\
\tilde{\lambda} + \tilde{m} &= 0
&&\text{in}\;\;\Omega\times\{1\},
\\
\D{H}\vect{\tilde{v}}
\defeq
\beta_v\D{B}[\vect{\tilde{v}}] + \D{L}[\vect{\tilde{b}}]
&= -\vect{g}
&&\text{in}\;\;\Omega,
\label{e:inc-control-eq-elim}
\end{align}
\end{subequations}

\noindent with periodic boundary conditions on $\p\Omega$. Here, $\vect{g}$ in~\eqref{e:inc-control-eq-elim} corresponds to the reduced gradient in~\eqref{e:control-eq-elim}. The operator $\D{B}$ is the second variation of the regularization model (see \secref{s:optimality-conditions}). The projection operator $\D{L}$ coincides with $\D{K}$ in~\eqref{e:project-body-force} if we consider the seminorms in~\eqref{e:reg-v} as a regularization operator. If we consider~\eqref{e:shear-control-reg-v} instead, we have
\begin{align}
\label{e:project-inc-body-force-nls}
\D{L}(\vect{v})[\vect{\tilde{b}},\vect{\tilde{v}}]
&=
\igrad
\D{M}^{-1}
\ilap^{-1}
\idiv
\left(
\idiv
2\beta_v
\left(
\hat{\eta}[\vect{v}]
+
\eta[\vect{v}]
\D{Q}[\vect{v}]
\right)
\D{E}[\vect{\tilde{v}}]
-
\vect{\tilde{b}}
\right)
+
\vect{\tilde{b}},
\end{align}

\noindent where the operator $\D{M}$ is as defined above and the operator $\D{Q}$ is given in~\eqref{e:2nd-var-reg-v-nls}. Our algorithm only operates on these reduced systems. We discuss its implementation, and by that the scheme to ultimately solve \eqref{e:opt-prob-v}, next.

\section{Numerics}
\label{s:numerics}

Our numerical scheme was originally described in~\cite{Mang:2015a}. We normalize the intensities of the images to $[0,1]$.  We use the trapezoidal rule for numerical quadrature and an explicit second order Runge--Kutta method for the numerical time integration of the hyperbolic PDEs in~\eqref{e:first-order-opt-elim} and~\eqref{e:second-order-elim}, respectively. Due to the conditional stability (CFL condition) we have to restrict the time step size $h_t$. Given that we invert for a stationary velocity field $\vect{v}$, we can modify the number of time steps $n_t$ as required. We use a pseudospectral discretization with a Fourier basis in space. This allows us to efficiently construct the inverse differential operators that appear in our formulation in \secref{s:elimination}.

Images are in general functions of bounded variation; our scheme cannot handle this type of discontinuity in the data. Accordingly, we assume that the images are adequately smooth. We ensure this numerically by pre-smoothing the data, a common strategy considered in many registration packages.\footnote{An inadequate smoothness of the data can deteriorate the accuracy of the solver. We numerically ensure stability by enforcing sufficient regularity of the data and the velocity field $\vect{v}$ based on a combination of pre-smoothing of the input data and an adequate choice of the regularization weight $\beta_v$. The adjoint equation plays a critical role; in our formulation, we have to differentiate the Lagrange multiplier. As pointed out in~\cite{Vialard:2012a}, we can ensure numerical feasibility if we ensure that the input images are adequately smooth. Further strategies include the use of another scheme for numerical time integration or to use a map-based formulation~\cite{Hart:2009a}. We will investigate this in the future.} We use a globalized, inexact~\cite{Dembo:1983a, Eisenstat:1996a}, preconditioned, matrix-free, reduced space (Gauss--)Newton--Krylov method for numerical optimization~\cite{Mang:2015a}. This scheme amounts to a sequential solution of the systems~\eqref{e:first-order-opt-elim} and~\eqref{e:second-order-elim}. The Newton step is in a general format given by
\begin{equation}
\label{e:reduced-space-kkt-system}
\D{H}^h_k \vect{\tilde{v}}^h_k = -\vect{g}^h_k,
\quad
\vect{v}^h_{k+1} = \vect{v}^h_k + \alpha_k\vect{\tilde{v}}_k^h,
\quad
k = 1,2,\ldots
\end{equation}

\noindent where $\D{H}^h_k \in \ns{R}^{n \times n}$, $n\in\ns{N}$, is a discrete representation of the reduced Hessian in~\eqref{e:inc-control-eq-elim} acting on the incremental control variable $\vect{\tilde{v}}^h_k\in\ns{R}^n$ at (outer) iteration $k$. The scheme is globalized via a backtracking line search subject to the Armijo condition (we use default parameters; see \cite[algorithm~3.1,\,p.~37]{Nocedal:2006a}). We iteratively solve~\eqref{e:reduced-space-kkt-system} using a PCG method. We refer to the solution of~\eqref{e:reduced-space-kkt-system} as \iname{inner iterations} (as opposed to the steps for updating $\vect{v}^h_k$, to which we refer to as \iname{outer iterations}). We ensure that the reduced space Hessian operator is positive definite by exploiting a Gauss--Newton approximation to the true Hessian. This corresponds to setting $\lambda$ in~\eqref{e:second-order-elim} to zero (see also~\cite{Mang:2015a}). The preconditioner for the reduced space Hessian is the inverse of the second variation of the regularization operator. This preconditioner has vanishing construction costs, due to the pseudospectral discretization in space.

We provide more details on this algorithm in the appendix \secref{s:algorithm}; we also refer to~\cite{Mang:2015a} for a detailed algorithmic study of our Newton--Krylov scheme in the context of large deformation diffeomorphic image registration; this includes a comparison to a preconditioned gradient descent scheme.

\section{Numerical Experiments}
\label{s:numerical-experiments}

We study the performance of the proposed formulation in different application scenarios, accounting for synthetic and real-world registration problems. All results reported in this study are limited to the two-dimensional case. Nothing in our formulation is specific to $d=2$; a three-dimensional implementation is ongoing work.

 We limit the first part of this study in \secref{s:res-quadratic-reg-v} to the quadratic regularization norms in~\eqref{e:reg-v}. Results for the nonlinear regularization model in~\eqref{e:shear-control-reg-v} are reported in \secref{s:res-shear-control}.

\begin{figure}
\centering
\includegraphics[width=0.99\textwidth]
{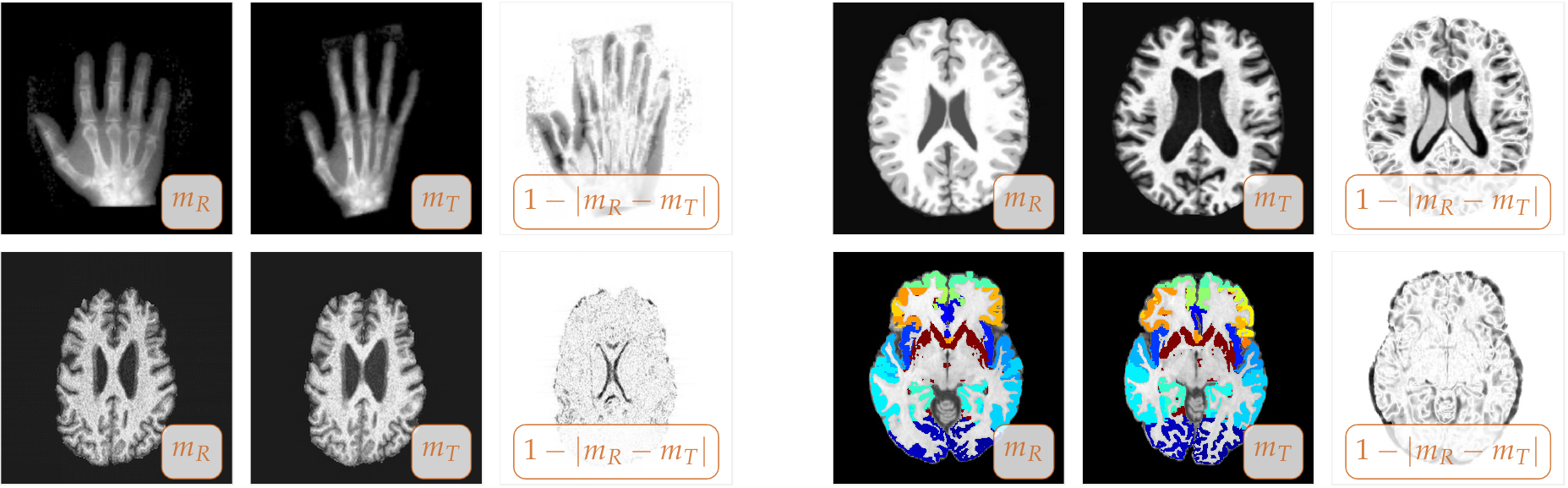}
\caption
{
Real-world two-dimensional registration problems. We display (from left to right) the reference image $m_R$ (fixed image), the template image $m_T$ (image to be registered), and a map of the residual differences between $m_R$ and $m_T$ before registration (for each set of images as indicated by the inset). Top left: benchmark registration problem \cite{Amit:1994a, Modersitzki:2004a, Modersitzki:2009a}; top right: intersubject registration problem; bottom left: longitudinal (intrasubject) registration problem; bottom right: intersubject registration problem~\cite{Christensen:2006a} (for the latter data we have a ground truth based on annotations: segmentations of 32 anatomical gray matter regions of interest; we overlay the associated label maps onto the reference and template image).
\label{f:regprob-image-collection-smooth}
}
\end{figure}

\subsection{Quadratic Regularization}
\label{s:res-quadratic-reg-v}

We report different measures of registration performance, with the aim to assess both, the fidelity of the registration as well as properties of the computed deformation map. We report results for different two-dimensional real-world registration problems.

\subsubsection{Registration Performance}

\ipoint{Purpose} We study registration quality as a function of the regularization parameters $\beta_v$ (smoothness) and $\beta_w$ (incompressiblity).

\ipoint{Setup} All images are registered in full resolution. No grid, scale, or parameter continuation is performed. We terminate the optimization if the relative reduction of the gradient is at least three orders of magnitude. We consider three two-dimensional, real-world registration problems: a benchmark problem based on \iname{hand images}\footnote{The images are taken from~\cite{Modersitzki:2009a}.}~\cite{Amit:1994a, Modersitzki:2004a, Modersitzki:2009a} as well as a multisubject and a serial\footnote{The data is available at \href{http://central.xnat.org}{\tt http://central.xnat.org} (\iname{open access series of imaging studies} (OASIS) longitudinal study; dataset 70; time point one and time point four). The images have been affinely preregistered~\cite{Klein:2010a} and skull stripped~\cite{Smith:2002a}.} (longitudinal) brain image registration problem (\iname{multisubject brain images} and \iname{serial brain images}). The initial images are displayed in~\figref{f:regprob-image-collection-smooth}.

The images have a grid size of $256\times256$. The number of time points is adapted as required by monitoring the CFL condition (initialized with $n_t = 2\max(\vect{n}_x)$). We vary $\beta_v$ and $\beta_w$ in steps of one order of magnitude ranging from $\num{1E-5}$ to $\num{1E-1}$, respectively. If we further reduce the regularization parameters the problem becomes computationally prohibitive (due to ill-conditioning) and numerically unstable (we not only violate the theoretical smoothness assumptions~\cite{Chen:2011a, Dupuis:1998a} but also approach regimes that are numerically unstable; this will eventually result in irregular, nondiffeomorphic mappings); smaller values for the regularization parameters require finer grids to resolve the problem. We terminate the optimization if the relative change of the $\ell^{\infty}$-norm of the reduced gradient is at least three orders of magnitude. We compare the designed framework for near-incompressible registration to plain smoothness regularization and a model for incompressible diffeomorphisms.

\begin{figure}
\centering
\includegraphics[width=0.8\textwidth]
{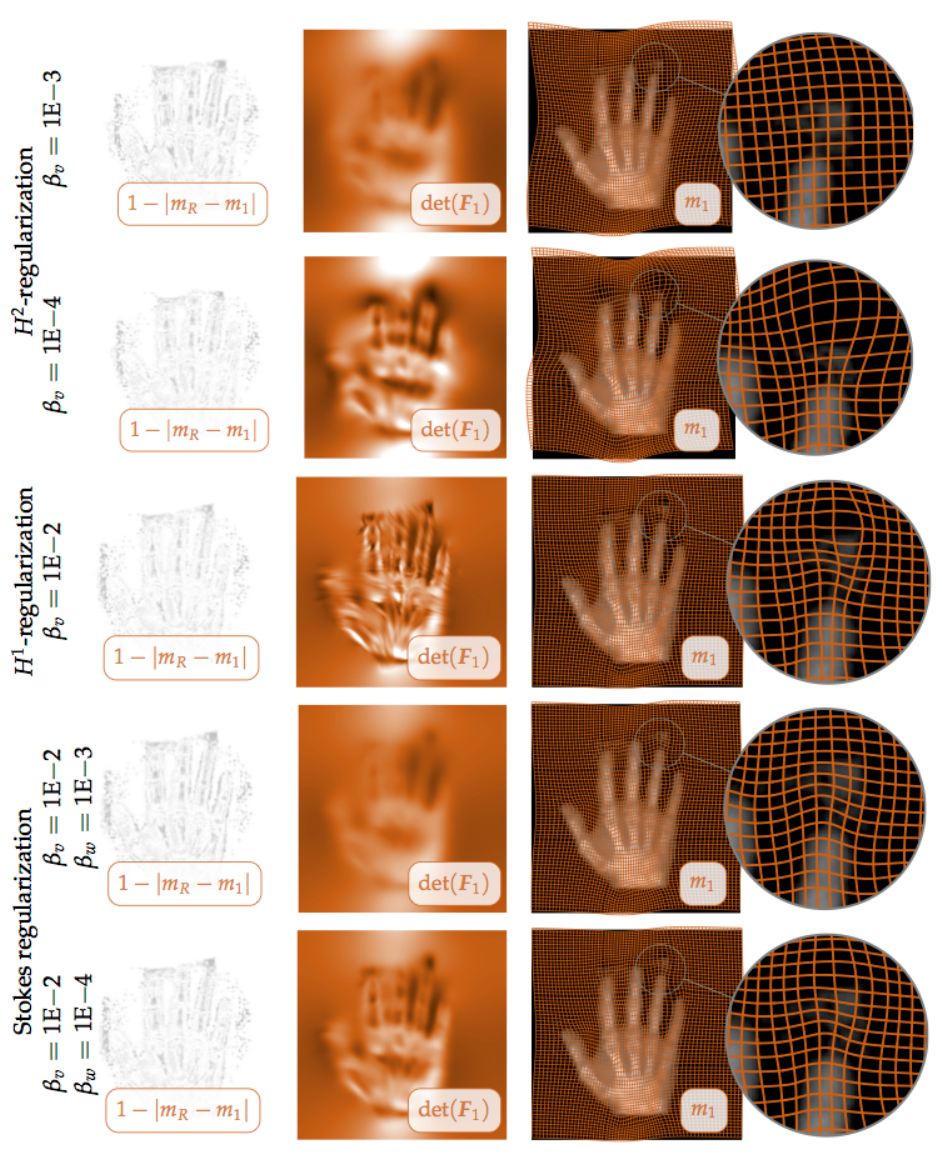}
\caption
{
Exemplary registration results for the \iname{hand images} (see~\figref{f:regprob-image-collection-smooth}). We report representative results from \tabref{t:quantitative-results-quadratic-regularization}. The first three rows show results for plain smoothness regularization ($\gamma=0$; first and second row: $H^2$-regularization; third row: $H^1$-regularization) for different choices of $\beta_v$ (top row: $\beta_v = \num{1E-3}$; second row: $\beta_v = \num{1E-4}$; third row: $\beta_v=\num{1E-2}$). The two rows from the bottom show results for a model with local adaptive compression ($H^1$-regularization; $\gamma=1$) for $\beta_v=\num{1E-2}$ and different choices for $\beta_w$ (bottom row: $\beta_2=\num{1E-4}$; second row from the bottom: $\beta_w=\num{1E-3}$). We show (from left to right) ($i$) the residual differences after registration, ($ii$) a map of the determinant of the deformation gradient (the values are reported in \tabref{t:quantitative-results-quadratic-regularization}; the color map is explained in \secref{s:visualization} of the appendix), ($iii$) the deformed template image $m_1$ with a grid in overlay, and ($iv$) a close up of the latter for a particular area of interest (as identified by the inset in the images).
}
\label{f:results-hands-exemplary}
\end{figure}

\begin{figure}
\centering
\includegraphics[width=0.8\textwidth]
{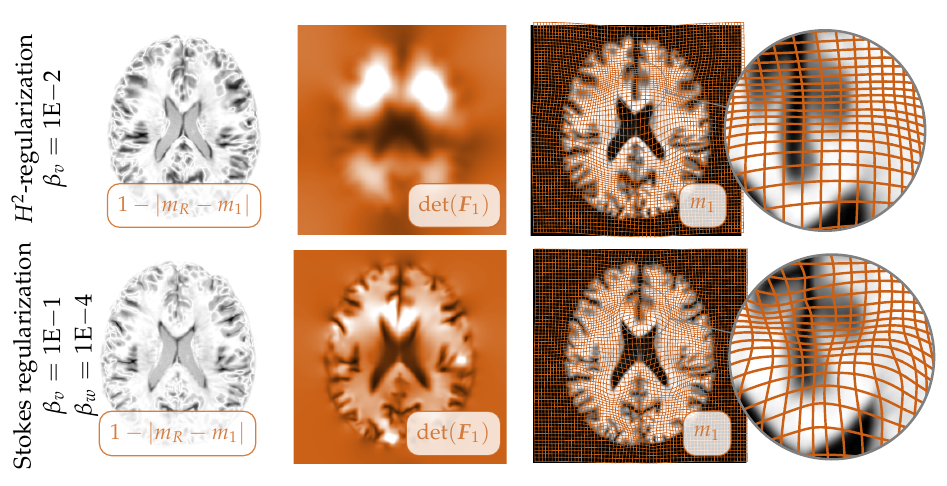}
\caption
{
Exemplary registration results for the \iname{multisubject brain images} (see \figref{f:regprob-image-collection-smooth}). We report representative results from~\tabref{t:quantitative-results-quadratic-regularization}. We report results for plain smoothness regularization (top row; $H^2$-regularization; $\gamma=0$; $\beta_v=\num{1E-2}$) and for a model with local adaptive compression (bottom row; $H^1$-regularization; $\gamma=1$; $\beta_v=\num{1E-1}$; $\beta_w=\num{1E-4}$). We display (from left to right) ($i$) a map of the residual differences after registration, ($ii$) a map of the determinant of the deformation gradient (the values are reported in \tabref{t:quantitative-results-quadratic-regularization}; information about the color map can be found in \secref{s:visualization} of the appendix), ($iii$) the deformed template image $m_1$ with a grid in overlay (to illustrate the deformation map $\vect{y}$), and ($iv$) a close up of the latter for a particular area of interest (as identified by the inset in the images).
}
\label{f:results-multi-subject-brains-exemplary}
\end{figure}

\begin{figure}
\centering
\includegraphics[width=0.8\textwidth]
{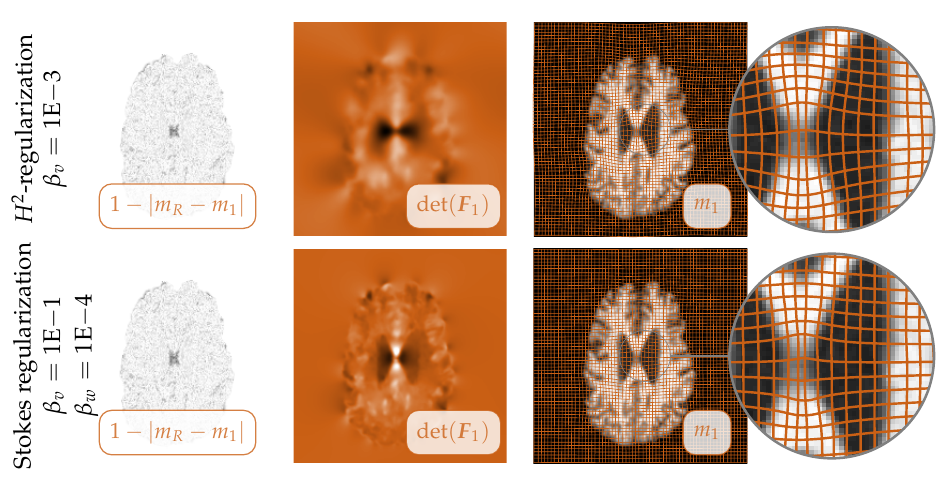}
\caption
{
Exemplary registration results for the \iname{serial brain images} (see \figref{f:regprob-image-collection-smooth}). We compare plain smoothness regularization based on an $H^2$-seminorm (top row; $\beta_v=\num{1E-3}$; $\gamma=0$) to the designed model with local adaptive compression (bottom row; $H^1$-regularization; $\beta_v=\num{1E-1}$; $\beta_w=\num{1E-4}$; $\gamma=1$). We display (from left to right) ($i$) a map of the residual differences between the reference image $m_R$ and the deformed template $m_1$, ($ii$) a map of the determinant of the deformation gradient (the values are reported in \tabref{t:quantitative-results-quadratic-regularization}; information about the color map can be found in \secref{s:visualization} of the appendix; notice that we changed the window to $[0.5,1.5]$, since the volume changes between the images are subtle), ($iii$) a deformed grid overlaid onto the deformed template image $m_1$ (to illustrate the deformation map $\vect{y}$) and ($iv$) a close up of the latter for a particular area of interest (as identified by the inset in the images).
}
\label{f:results-serial-brains-exemplary}
\end{figure}

\ipoint{Results} Quantitative results are summarized in~\tabref{t:quantitative-results-quadratic-regularization}. Exemplary results for the \iname{hand images}, the \iname{multisubject brain images} and the \iname{serial brain images} are illustrated in~\figref{f:results-hands-exemplary}, \figref{f:results-multi-subject-brains-exemplary} and \figref{f:results-serial-brains-exemplary}, respectively.

\setcounter{runidnum}{0}
\begin{table}
\caption
{
Quantitative analysis of registration performance as a function of the regularization parameters $\beta_v$ and $\beta_w$. The registration problem are the \iname{hand images}, the \iname{multisubject brain images}, and the \iname{serial brain images} (see \figref{f:regprob-image-collection-smooth}). We report results for different regularization schemes: We report ($i$) results for smoothness regularization without a constraint on the divergence of the velocity field ($H^1$- and $H^2$-regularization on $\vect{v}$; $\gamma=0$), ($ii$) for incompressible diffeomorphisms ($H^1$-regularization on $\vect{v}$; $\gamma=1$), and ($iii$) the proposed model with local adaptive compression ($H^1$-regularization on $\vect{v}$ and $w$; $\gamma=1$). We report values for ($i$) the number of Hessian matrix vector products ($n_{\text{matvec}}$), ($ii$) the number of hyperbolic PDE solves ($n_{\text{PDE}}$), ($iii$) the relative reduction of the gradient ($\|\vect{g}^\star\|_{\rm rel}$), ($iv$) the relative reduction of the mismatch ($\|\vect{r}^\star\|_{\rm rel}$), ($v$) $\min$, $\max$ and $\operatorname{mean}$ values of the determinant of the deformation gradient $J$, and ($vi$) the $\max$ and $\operatorname{mean}$ distance of the deformation gradient from identity $D$ (indicating the distance from a completely rigid mapping). The definitions of these measures can be found in \tabref{t:registration-performance-measures} in \secref{s:measures-reg-performance} in the appendix.
\label{t:quantitative-results-quadratic-regularization}
}
\centering
\begin{scriptsize}
\begin{tabular}{rrrrrgrrggrrrr}
\toprule
  run
& $\|\vect{v}\|^2_{\fs{V}}$
& $\gamma$
& $\beta_v$
& $\beta_w$
& $n_{\text{mv}}$
& $n_{\text{PDE}}$
& $\|\vect{g}^{\star}\|_{\text{rel}}$
& $\|\vect{r}^{\star}\|_{\text{rel}}$
& $\min(J)$
& $\operatorname{mean}(J)$
& $\max(J)$
& $\operatorname{mean}(D)$
& $\max(D)$
\\\midrule
\multicolumn{14}{c}{\it hand images}
\\\midrule
  \runid
& $H^2$
& 0
& \num[round-precision=1]{1.000000e-01}
& n/a
& 145
& 328
& \num{6.507795e-04}
& \num{2.572041e-01}
& \num{7.764302e-01}
& \num{1.019638e+00}
& \num{1.334697e+00}
& \num{1.592040e-01}
& \num{3.366220e-01}
\\
  \runid
&
&
& \num[round-precision=1]{1.000000e-02}
&
& 475
& 999
& \num{7.950515e-04}
& \num{1.224387e-01}
& \num{6.901432e-01}
& \num{1.047543e+00}
& \num{1.726413e+00}
& \num{2.527154e-01}
& \num{5.535723e-01}
\\
  \runid
&
&
& \num[round-precision=1]{1.000000e-03}
&
& 2435
& 4972
& \num{8.856170e-04}
& \num{8.406727e-02}
& \num{5.984229e-01}
& \num{1.068857e+00}
& \num{1.986559e+00}
& \num{3.173773e-01}
& \num{8.794028e-01}
\\
  \runid
&
&
& \num[round-precision=1]{1.000000e-04}
&
& 15189
& 30648
& \num{9.710304e-04}
& \num{5.456891e-02}
& \num{2.895222e-01}
& \num{1.107307e+00}
& \num{2.448447e+00}
& \num{4.074396e-01}
& \num{1.319238e+00}
\\\midrule
  \runid
& $H^1$
& 0
& \num[round-precision=1]{1.000000e-01}
& n/a
& 187
& 419
& \num{7.949208e-04}
& \num{1.599736e-01}
& \num{6.538843e-01}
& \num{1.019852e+00}
& \num{1.638156e+00}
& \num{1.607336e-01}
& \num{5.108978e-01}
\\
  \runid
&
&
& \num[round-precision=1]{1.000000e-02}
&
& 1560
& 3291
& \num{9.461070e-04}
& \num{5.562520e-02}
& \num{2.716105e-01}
& \num{1.057555e+00}
& \num{2.607476e+00}
& \num{2.634732e-01}
& \num{2.142862e+00}
\\\midrule
  \runid
& $H^1$
& 1
& \num[round-precision=1]{1.000000e-01}
& n/a
& 392
& 868
& \num{9.283519e-04}
& \num{3.888012e-01}
& \num{9.999540e-01}
& \num{1.000000e+00}
& \num{1.000046e+00}
& \num{1.639276e-01}
& \num{7.733198e-01}
\\
  \runid
&
&
& \num[round-precision=1]{1.000000e-02}
&
& 1011
& 2116
& \num{8.737323e-04}
& \num{2.379146e-01}
& \num{9.985685e-01}
& \num{1.000000e+00}
& \num{1.001418e+00}
& \num{2.985403e-01}
& \num{1.556337e+00}
\\\midrule
  \runid
& $H^1$
& 1
& \num[round-precision=1]{1.000000e-01}
& \num[round-precision=1]{1.000000e-01}
& 306
& 675
& \num{8.982489e-04}
& \num{2.626872e-01}
& \num{8.999427e-01}
& \num{1.004499e+00}
& \num{1.141837e+00}
& \num{1.468633e-01}
& \num{5.783271e-01}
\\
  \runid
&
&
&
& \num[round-precision=1]{1.000000e-02}
& 273
& 603
& \num{7.748664e-04}
& \num{1.946754e-01}
& \num{7.996145e-01}
& \num{1.014083e+00}
& \num{1.337582e+00}
& \num{1.538733e-01}
& \num{4.511429e-01}
\\
  \runid
&
&
&
& \num[round-precision=1]{1.000000e-03}
& 241
& 537
& \num{8.804795e-04}
& \num{1.703212e-01}
& \num{7.309260e-01}
& \num{1.018279e+00}
& \num{1.513027e+00}
& \num{1.591725e-01}
& \num{4.546416e-01}
\\
  \runid
&
&
&
& \num[round-precision=1]{1.000000e-04}
& 212
& 472
& \num{7.891485e-04}
& \num{1.621329e-01}
& \num{6.998720e-01}
& \num{1.019506e+00}
& \num{1.614504e+00}
& \num{1.603560e-01}
& \num{4.810370e-01}
\\
  \runid
&
&
&
& \num[round-precision=1]{1.000000e-05}
& 192
& 429
& \num{8.488561e-04}
& \num{1.602730e-01}
& \num{6.612167e-01}
& \num{1.019801e+00}
& \num{1.635024e+00}
& \num{1.606649e-01}
& \num{5.043789e-01}
\\
  \runid
&
&
& \num[round-precision=1]{1.000000e-02}
& \num[round-precision=1]{1.000000e-01}
& 1819
& 3756
& \num{9.962345e-04}
& \num{1.537354e-01}
& \num{9.435955e-01}
& \num{1.003149e+00}
& \num{1.179269e+00}
& \num{2.811630e-01}
& \num{1.406140e+00}
\\
  \runid
&
&
&
& \num[round-precision=1]{1.000000e-02}
& 1466
& 3025
& \num{9.807377e-04}
& \num{1.029497e-01}
& \num{8.389104e-01}
& \num{1.016092e+00}
& \num{1.429109e+00}
& \num{2.644362e-01}
& \num{1.101264e+00}
\\
  \runid
&
&
&
& \num[round-precision=1]{1.000000e-03}
& 1385
& 2871
& \num{9.718678e-04}
& \num{8.289399e-02}
& \num{6.686035e-01}
& \num{1.034995e+00}
& \num{1.696425e+00}
& \num{2.514603e-01}
& \num{1.049177e+00}
\\
  \runid
&
&
&
& \num[round-precision=1]{1.000000e-04}
& 1483
& 3083
& \num{9.353795e-04}
& \num{6.737001e-02}
& \num{5.137283e-01}
& \num{1.046796e+00}
& \num{1.841633e+00}
& \num{2.545157e-01}
& \num{1.242277e+00}
\\
  \runid
&
&
&
& \num[round-precision=1]{1.000000e-05}
& 1404
& 2927
& \num{9.230329e-04}
& \num{5.865952e-02}
& \num{3.440553e-01}
& \num{1.053888e+00}
& \num{2.339445e+00}
& \num{2.594844e-01}
& \num{1.522920e+00}
\\\midrule
\multicolumn{14}{c}{\it multisubject brain images}
\\\midrule
  \runid
& $H^2$
& 0
& \num[round-precision=1]{1.000000e-01}
& n/a
& 7538
& 16223
& \num{9.883048e-04}
& \num{7.348915e-01}
& \num{5.916770e-01}
& \num{1.016491e+00}
& \num{1.328465e+00}
& \num{1.247653e-01}
& \num{3.655436e-01}
\\
  \runid
&
&
& \num[round-precision=1]{1.000000e-02}
&
& 5685
& 11855
& \num{9.727689e-04}
& \num{4.490202e-01}
& \num{3.076747e-01}
& \num{1.093236e+00}
& \num{2.678074e+00}
& \num{2.792158e-01}
& \num{1.283676e+00}
\\\midrule
  \runid
& $H^1$
& 1
& \num[round-precision=1]{1.000000e-01}
& \num[round-precision=1]{1.000000e-01}
& 3050
& 6722
& \num{9.990809e-04}
& \num{7.216611e-01}
& \num{6.918409e-01}
& \num{1.006264e+00}
& \num{1.176897e+00}
& \num{1.466081e-01}
& \num{6.748061e-01}
\\
  \runid
&
&
&
& \num[round-precision=1]{1.000000e-02}
& 3613
& 7971
& \num{9.840200e-04}
& \num{4.904388e-01}
& \num{4.179794e-01}
& \num{1.036812e+00}
& \num{1.844239e+00}
& \num{2.049549e-01}
& \num{1.109927e+00}
\\
  \runid
&
&
&
& \num[round-precision=1]{1.000000e-03}
& 2107
& 4627
& \num{9.726776e-04}
& \num{3.521920e-01}
& \num{3.055548e-01}
& \num{1.070993e+00}
& \num{2.145434e+00}
& \num{2.471565e-01}
& \num{1.710717e+00}
\\
  \runid
&
&
&
& \num[round-precision=1]{1.000000e-04}
& 2305
& 4985
& \num{9.792073e-04}
& \num{2.946096e-01}
& \num{1.673288e-01}
& \num{1.087696e+00}
& \num{2.514142e+00}
& \num{2.623297e-01}
& \num{1.856906e+00}
\\
  \runid
&
&
&
& \num[round-precision=1]{1.000000e-05}
& 2703
& 5858
& \num[round-precision=2]{9.999033e-04}
& \num{2.787460e-01}
& \num{1.257547e-01}
& \num{1.093080e+00}
& \num{2.717418e+00}
& \num{2.654181e-01}
& \num{1.874935e+00}
\\\midrule
\multicolumn{14}{c}{\it serial brain images}
\\\midrule
  \runid
& $H^2$
& 0
& \num[round-precision=1]{1.000000e-01}
& n/a
& 67
& 151
& \num{7.244139e-04}
& \num{6.720449e-01}
& \num{9.039711e-01}
& \num{1.000502e+00}
& \num{1.078836e+00}
& \num{2.065292e-02}
& \num{1.184992e-01}
\\
  \runid
&
&
& \num[round-precision=1]{1.000000e-02}
&
& 169
& 350
& \num{3.986599e-04}
& \num{4.587686e-01}
& \num{7.385633e-01}
& \num{1.002006e+00}
& \num{1.154629e+00}
& \num{3.940637e-02}
& \num{3.255374e-01}
\\
  \runid
&
&
& \num[round-precision=1]{1.000000e-03}
&
& 1286
& 2607
& \num{9.100350e-04}
& \num{2.979252e-01}
& \num{5.136548e-01}
& \num{1.005371e+00}
& \num{1.298239e+00}
& \num{7.117103e-02}
& \num{9.226124e-01}
\\
  \runid
&
&
& \num[round-precision=1]{1.000000e-04}
&
& 13772
& 27696
& \num{1.245593e-04}
& \num{1.873162e-01}
& \num{2.624056e-01}
& \num{1.014218e+00}
& \num{1.739953e+00}
& \num{1.259960e-01}
& \num{1.951390e+00}
\\\midrule
  \runid
& $H^1$
& 1
& \num[round-precision=1]{1.000000e-01}
& \num[round-precision=1]{1.000000e-01}
& 103
& 227
& \num{3.411192e-04}
& \num{5.445334e-01}
& \num{9.158836e-01}
& \num{1.000310e+00}
& \num{1.059657e+00}
& \num{3.689540e-02}
& \num{3.583418e-01}
\\
  \runid
&
&
&
& \num[round-precision=1]{1.000000e-02}
& 78
& 173
& \num{8.495168e-04}
& \num{4.180768e-01}
& \num{7.869497e-01}
& \num{1.001279e+00}
& \num{1.134772e+00}
& \num{3.821299e-02}
& \num{4.311000e-01}
\\
  \runid
&
&
&
& \num[round-precision=1]{1.000000e-03}
& 107
& 235
& \num{8.988925e-04}
& \num{3.422411e-01}
& \num{6.761942e-01}
& \num{1.002404e+00}
& \num{1.285009e+00}
& \num{4.003887e-02}
& \num{7.303273e-01}
\\
  \runid
&
&
&
& \num[round-precision=1]{1.000000e-04}
& 117
& 258
& \num{8.473846e-04}
& \num{3.083017e-01}
& \num{6.233797e-01}
& \num{1.003305e+00}
& \num{1.694508e+00}
& \num{4.159524e-02}
& \num{1.103124e+00}
\\
  \runid
&
&
&
& \num[round-precision=1]{1.000000e-05}
& 177
& 390
& \num{8.479657e-04}
& \num{2.977446e-01}
& \num{6.037333e-01}
& \num{1.003724e+00}
& \num{1.881677e+00}
& \num{4.256186e-02}
& \num{1.291362e+00}
\\\bottomrule
\end{tabular}
\end{scriptsize}
\end{table}

\ipoint{Observations} The most important observations are the following: Augmenting smoothness regularization with a constraint on $\idiv\vect{v}$ with a nonzero right-hand side $w$ (mass source) allows us to control the magnitude of the determinant of the deformation gradient without compromising registration quality. We avoid oversmoothing of the deformation map $\vect{y}$.

Enforcing incompressibility up to the numerical error is not adequate for the considered registration problems. This is also reflected by the residual differences reported in \tabref{t:quantitative-results-quadratic-regularization}. Using a plain $H^1$-seminorm as a regularization model (with no control on $\idiv\vect{v}$) can be delicate: small variations in the regularization parameter $\beta_v$ yield strong variations in the determinant of the deformation gradient. The divergence constraint allows us to better control the mapping. The trend of the values for $\det(\mat{F}_1)$ as a function of $\beta_v$ and $\beta_w$ demonstrates that we can precisely control the regularity properties of the mapping $\vect{y}$.

In some cases we can---as compared to plain smoothness regularization---significantly reduce the variations of the determinant of the deformation gradient without compromising registration quality. For instance in run \#16 in \tabref{t:quantitative-results-quadratic-regularization} we set $\beta_v$ to $\num{1E-2}$ and $\beta_w$ to $\num{1E-3}$ and obtain an $L^2$-distance of $\num{8.289399e-02}$ with $\det(\mat{F}_1)\in[\num{6.686035e-01},\num{1.696425e+00}]$. The maximum and mean distance of the deformation gradient from identity is $\num{2.514603e-01}$ and $\num{1.049177e+00}$, respectively. If we want to obtain a similar residual using plain smoothness regularization on the basis of an $H^1$-seminorm, we have to set $\beta_v$ to $\num{1E-2}$ (run \#6 in \tabref{t:quantitative-results-quadratic-regularization}). This results in a relative change of the $L^2$-distance of $\num{5.562520e-02}$. However, the variation of the determinant of the deformation gradient is larger with $\det(\mat{F}_1)\in[\num{2.716105e-01}, \num{2.607476e+00}]$. If we use an $H^2$-seminorm we achieve a similar mismatch (relative reduction of the $L^2$-distance by $\num{8.406727e-02}$) for $\beta_v=\num{1E-3}$ (run \#3 in \tabref{t:quantitative-results-quadratic-regularization}). The variation in the determinant of the deformation gradient is slightly larger with $\det(\mat{F}_1)\in[\num{5.984229e-01}, \num{1.986559e+00}]$ (as compared to the near-incompressible case).

Careful visual inspection of the results in \figref{f:results-hands-exemplary} confirms these findings. The residual differences are very similar for all models. We can also see that if we set $\beta_v$ to $\num{1E-2}$ or $\num{1E-4}$ for plain $H^1$- and $H^2$-regularization (i.e.\ without additional constraint on the divergence of $\vect{v}$; runs \#4 and \#6 in \tabref{t:quantitative-results-quadratic-regularization}), we seem to overfit the data; the mapping becomes more and more ill-behaved. By setting $\beta_v$ to $\num{1E-3}$ for the $H^2$-seminorm or by using a near-incompressible model with $\beta_v=\num{1E-2}$ a nice diffeomorphism is obtained.

For the \iname{hand images} we obtain an equivalent performance for the $H^1$- and $H^2$-regularization models since the mapping between both images is rather smooth. This is different for the \iname{multisubject brain images} (see \figref{f:results-multi-subject-brains-exemplary}). The $H^2$-seminorm yields a nice diffeomorphic map but $\vect{y}$ also appears to be overly smooth. That is, we observe a strong blurring in the map of the determinant of the deformation gradient. Thus, it is not possible to recover fine features in the deformation field. The same behavior can be observed for the \iname{serial brain images}. Although the residual differences are very similar for the $H^1$- and the $H^2$-seminorm, we obtain mappings that are locally very different. If we use an $H^1$-regularization we can recover much more localized features in the deformation map. These local changes could be of interest in a subsequent analysis of the local deformation properties (volume changes; deformation based morphometry). Also note that the mean values for the determinant of the deformation gradient are closer to one (i.e.\ volume is more likely preserved) as compared to plain smoothness regularization.

When switching from $H^1$- to an $H^2$-regularization model we have to reduce $\beta_v$ by one order of magnitude to obtain a similar mismatch. Note that the computational complexity of our scheme is currently not mesh independent; the rate of convergence deteriorates significantly if we reduce $\beta_v$ as judged by the number of Hessian matrix vector products and hyperbolic PDE solves.

\ipoint{Conclusions} Using an $H^2$-regularization model---a common choice in large deformation diffeomorphic registration algorithms---yields well behaved mappings. However, we might loose local features (fine structures) in the deformation map due to oversmoothing. These features could be of importance for a subsequent analysis of the deformation map. Empirically, we observed that we have to reduce $\beta_v$ by one order of magnitude for the $H^2$-regularization model as compared to the $H^1$-seminorm to obtain a similar mismatch. If we reduce $\beta_v$ significantly the computational work load for the inversion (with the defined tolerance) becomes prohibitive. If we switch to an $H^1$-seminorm as a regularization model we can resolve fine features in the deformation. Introducing a constraint on the divergence of the velocity field with a nonzero mass source $w$ allows us to explicitly control the magnitude of the determinant of the deformation gradient without compromising registration quality. This relaxation of the incompressibility constraint is critical to make our model applicable across a wide range of registration problems.

\subsubsection{Validation and Comparison}

\ipoint{Purpose} We compare and validate registration performance of our algorithm. We report results for our formulation and the diffeomorphic {\tt DEMONS} algorithm~\cite{Vercauteren:2007a,Vercauteren:2009a}.

\ipoint{Setup} Our evaluation is based on the data of the \iname{Nonrigid Image Registration Evaluation Project} ({\tt NIREP}) \cite{Christensen:2006a}.\footnote{The data is available at \href{http://nirep.org}{\tt http://nirep.org} and described in detail in \cite{Christensen:2006a}.} {\tt NIREP} is a standardized data repository for the validation of deformable image registration algorithms; we refer to~\cite{Christensen:2006a} for details. We use the datasets {\tt na01} and {\tt na02} to study registration performance as a function of regularization parameters and norms (see \figref{f:regprob-image-collection-smooth}). Since our implementation is only a two-dimensional prototype, we extract a single slice from both volumes (axial slice 128) and resample the data to a resolution of $\vect{n}_x = (128,150)^\T$ (using a spline interpolation model for the image data and a nearest neighbor interpolation model for the label maps). For simplicity, we do not report results for the 32 individual labels but combine them to a single gray matter label map. We report the JSC, the DSC, the FPE, and the FNE (see \tabref{t:registration-performance-measures} in \secref{s:measures-reg-performance} for the definitions). We also report values for the determinant of the deformation gradient. We limit the evaluation of the deformation gradient to the area occupied by the brain (identified by thresholding).

We perform a parameter continuation in $\beta_v$ with bounds on the minimal tolerable determinant of the deformation gradient (binary search; see \cite{Mang:2015a} for details) for our algorithm starting with $\beta_v=\num{1E-0}$. We perform this binary search for fixed values of $\beta_w$ (we vary $\beta_w$ by one order of magnitude, starting with $\beta_w=\num{1E-1}$; i.e., we perform an exhaustive search for the second regularization parameter. We do not perform an additional grid or scale continuation. We terminate the optimization if the relative reduction in the gradient is two orders of magnitude or more.\footnote{We use a larger tolerance in these experiments because we did not observe any differences in the results if we turned to smaller tolerances.} Once we have obtained the velocity field (i.e., we solved the inverse problem with an estimated, optimal combination of regularization parameters $\beta_v$ and $\beta_w$) we compute the determinant of the deformation gradient and transport the label maps (as a post-processing step). We do so at a grid size of $4\vect{n}_x$ to be able to fully resolve the problem. Before we solve the transport problem we smooth the label maps using a Gaussian filter with standard deviation of $3\vect{h}_x$ (to avoid Gibbs phenomena). We threshold the transported label maps at a threshold of 0.5 to obtain binary labels and map them to the original resolution level by injection. We subsequently compute the overlap between the reference and transported template label maps.

The publicly available {\tt DEMONS} algorithm\footnote{We use the public implementation found at \href{http://hdl.handle.net/1926/510}{http://hdl.handle.net/1926/510} (see \cite{Vercauteren:2007a}; compiled with ITK4).} does not provide any stopping conditions other than the number of iterations. We tested several settings for the number of iterations in combination with a varying number of multiresolution levels. We observed that an increase in the number of iterations does not necessarily improve the obtained results; as a matter of fact the results can deteriorate for certain regularization parameter combinations if the iterations per resolution level are increased (i.e., the algorithm diverges). After these initial experiments we decided to fix the number of iterations and any other settings to the default values suggested in the documentation of the code (three grid resolution levels with 15, 10, and 5 iterations, respectively, with a diffeomorphic update rule based on an exponential map~\cite{Vercauteren:2009a} and symmetric gradient forces). We study the registration accuracy of the {\tt DEMONS} algorithm as a function of the two regularization parameters $\sigma_u$ (smoothing for the update field; \iname{fluid-like regularization}) and $\sigma_d$ (smoothing for the deformation field; \iname{diffusive regularization}). We report results for \bipa\item the mixed case, \item pure \iname{diffusive regularization}, and \item pure \iname{fluid-like regularization}\eipa. We report results for high data fidelity and low deformation regularity as well as results for low data fidelity and high deformation regularity. We apply the obtained mapping to the label maps using a nearest neighbor interpolation model.

For our algorithm, we use the parameters from the run above to extend our analysis to the remaining {\tt NIREP} datasets. That is, we use the values for the regularization parameters $\beta_v$ and $\beta_w$ that resulted in the best DSC scores for the registration between {\tt na01} and {\tt na02} to register all the remaining images to {\tt na01}. For the {\tt   DEMONS} algorithm we use two parameter settings for $\sigma_d$ and $\sigma_u$: one setting that results in the best DSC score (under the constraint that the deformation has to be diffeomorphic) and one setting that matches the determinant of the deformation gradient delivered by our method.

\begin{figure}
\centering
\includegraphics[width=0.8\textwidth]
{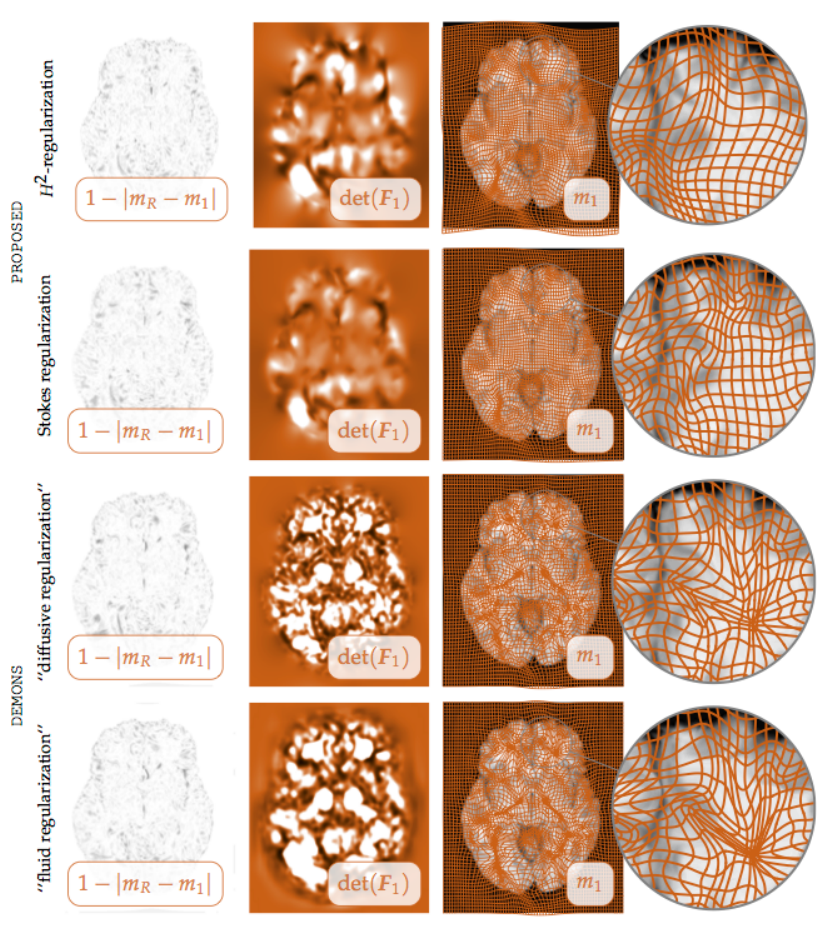}
\caption
{
Performance evaluation. We report qualitative results for the {\tt NIREP} data sets for our method (top row and middle row) and the {\tt DEMONS} algorithm (bottom row). The displayed results correspond to the runs reported in \tabref{t:performance-evaluation-nirep-parametersweep} ($H^2$-regularization: \runref{26}; linear Stokes regularization: \runref{31}; {\tt DEMONS}: \runref{5} (diffusive regularization) and \runref{22} (fluid-like regularization); we use run \runref{26} instead of the best run (\runref{27}) for the $H^2$-regularization because it has a similar DSC as we obtain for the linear Stokes case). We report for each run (from left to right) ($i$) the residual differences after registration, ($ii$) a map of the determinant of the deformation gradient, and ($iii$) a deformed grid overlaid onto the deformed template image. We computed the displayed results on a finer grid ($512\times600$) as compared to the one we have used to solve the optimization problem ($128\times150$) to be able to visualize the obtained deformation maps accurately (since MATLAB uses linear interpolation for visualization and we are using a spectral basis).
}
\label{f:results-nirep-exemplary}
\end{figure}

\begin{figure}
\centering
\includegraphics[width=0.95\textwidth]
{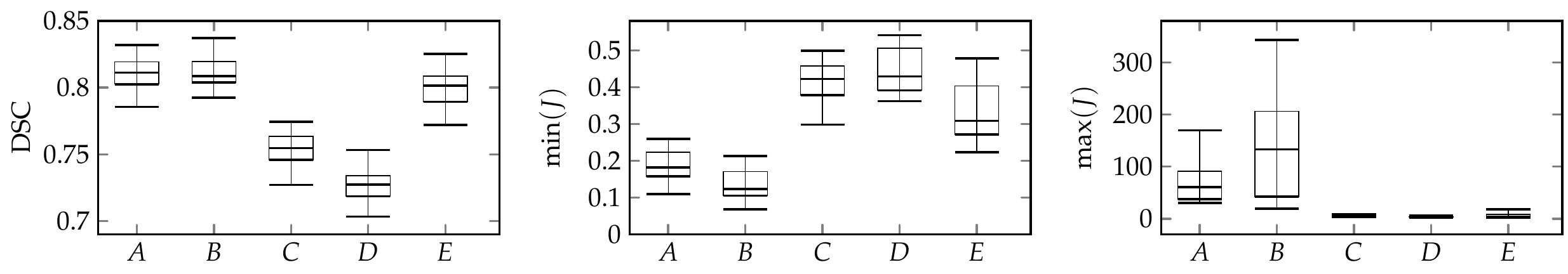}
\caption
{
Statistical quantification of registration quality across multiple datasets. We show box whisker plots to summarize the results reported in \tabref{t:performance-evaluation-nirep}. We compare registration quality for the {\tt DEMONS} algorithm and the {\tt PROPOSED} algorithm.  We report results for $A$: {\tt DEMONS}, $\sigma_u = 0$, $\sigma_d = 1$; $B$: {\tt DEMONS}, $\sigma_u = 2$, $\sigma_d = 0$; $C$: {\tt DEMONS}, $\sigma_u = 0$, $\sigma_d = 2$; $D$: {\tt DEMONS}, $\sigma_u = 4$, $\sigma_d = 0$; $E$: {\tt PROPOSED}, linear Stokes regularization, $\beta_v=\num{4.937500000000E-03}$, $\beta_w=\num{1E-4}$. These parameters are the same as in~\tabref{t:performance-evaluation-nirep}; these parameters deliver results that are consistent with our formulation in terms of the DSC scores or the values for the determinant of the deformation gradient. We report the DSC scores (left) and the smallest (middle) and largest (right) values for the determinant of the deformation gradient.
\label{f:results-nirep-quantitative-evaluation-boxplots}
}
\end{figure}

\ipoint{Results} We summarize the quantitative results for the datasets {\tt na01} and {\tt na02} in \tabref{t:performance-evaluation-nirep-parametersweep}. We provide qualitative results in \figref{f:results-nirep-exemplary}. We report the results for the remaining NIREP datasets in \tabref{t:performance-evaluation-nirep}; we summarize this experiment in \figref{f:results-nirep-quantitative-evaluation-boxplots}.

\setcounter{runidnum}{0}
\begin{table}
\caption
{
Performance evaluation. We compare registration quality between the diffeomorphic {\tt DEMONS} algorithm (top block) and the proposed algorithm (bottom block; plain $H^2$-smoothness regularization and linear Stokes regularization (LS)). We perform a parameter continuation in $\beta_v$ with lower bounds on the determinant of the deformation gradient for our approaches (the bound is 0.1 for LS and 0.5, 0.2, and 0.1 for the $H^2$-regularization model). For the {\tt DEMONS} algorithm we study registration quality as a function of the regularization parameters $\sigma_u$ (standard deviation for the smoothing of the update field) and $\sigma_d$ (standard deviation for the smoothing of the deformation field); i.e., we perform an exhaustive parameter search. We report values for the JSC, the DSC, the FPE, and the FNE. We also report $\min$, $\operatorname{mean}$, and $\max$ values for the determinant of the deformation gradient $J$.  We highlight the best results for each approach in bold. We highlight results that are identified to be nondiffeomorphic in faint gray color.
}
\label{t:performance-evaluation-nirep-parametersweep}
\centering
\begin{scriptsize}
\begin{tabular}{rrrrrgrrgrg}
\toprule
  run
&
& $\beta_v; \sigma_u$
& $\beta_w; \sigma_d$
& JSC
& DSC
& FPE
& FNE
& $\min(J)$
& $\operatorname{mean}(J)$
& $\max(J)$
\\\midrule
\multicolumn{11}{c}{\tt DEMONS}
\\\midrule
  \runid
& EXP
& \num{0.000000000000E+00}
& \num{4.000000000000E+00}
& \num{5.231686214190E-01}
& \num{6.869477404696E-01}
& \num{3.502076843198E-01}
& \num{2.936137071651E-01}
& \num{7.416120731755E-01}
& \num{1.077564819452E+00}
& \num{1.547330248312E+00}
\\
  \runid
&
& \num{0.000000000000E+00}
& \num{3.000000000000E+00}
& \num{5.575168583895E-01}
& \num{7.159047497772E-01}
& \num{3.089304257529E-01}
& \num{2.702492211838E-01}
& \num{6.179999788122E-01}
& \num{1.096551488689E+00}
& \num{2.287806442105E+00}
\\
  \runid
&
& \num{0.000000000000E+00}
& \num{2.000000000000E+00}
& \num{6.004998958550E-01}
& \num{7.503904216554E-01}
& \num{2.463655244029E-01}
& \num{2.515576323988E-01}
& \num{4.224464519978E-01}
& \num{1.186697645856E+00}
& \num{7.123601060829E+00}
\\
  \runid
&
& \num{0.000000000000E+00}
& \num{1.500000000000E+00}
& \num{6.378156701921E-01}
& \num{7.788613600422E-01}
& \num{2.027518172378E-01}
& \num{2.328660436137E-01}
& \num{2.972412656043E-01}
& \num{1.401283973047E+00}
& \num{2.362688763334E+01}
\\
 \bfseries\runid
&
&\bfseries\num{0.000000000000E+00}
&\bfseries\num{1.000000000000E+00}
&\bfseries\num{6.647501105705E-01}
&\bfseries\num{7.986184909671E-01}
&\bfseries\num{1.739356178609E-01}
&\bfseries\num{2.196261682243E-01}
&\bfseries\num{2.174175403328E-01}
&\bfseries\num{1.826167477435E+00}
&\bfseries\num{6.063498729499E+01}
\\
  \color{gray}\runid
&
& \color{gray}\num{0.000000000000E+00}
& \color{gray}\num{5.000000000000E-01}
& \color{gray}\num{7.334878331402E-01}
& \color{gray}\num{8.462566844920E-01}
& \color{gray}\num{1.201973001038E-01}
& \color{gray}\num{1.783489096573E-01}
& \color{gray}\num{-4.580603918460E+01}
& \color{gray}\num{3.755706151176E+01}
& \color{gray}\num{2.157974013809E+04}
\\
  \runid
&
& \num{5.000000000000E-01}
& \num{4.000000000000E+00}
& \num{5.225831251201E-01}
& \num{6.864428174703E-01}
& \num{3.507268951194E-01}
& \num{2.941329179647E-01}
& \num{7.421528928279E-01}
& \num{1.077466989543E+00}
& \num{1.549974348732E+00}
\\
  \runid
&
& \num{5.000000000000E-01}
& \num{3.000000000000E+00}
& \num{5.561719833564E-01}
& \num{7.147950089127E-01}
& \num{3.102284527518E-01}
& \num{2.712876427830E-01}
& \num{6.178134722172E-01}
& \num{1.096420844805E+00}
& \num{2.279960072790E+00}
\\
  \runid
&
& \num{5.000000000000E-01}
& \num{2.000000000000E+00}
& \num{5.978373882304E-01}
& \num{7.483081728267E-01}
& \num{2.484423676012E-01}
& \num{2.536344755971E-01}
& \num{4.222185905843E-01}
& \num{1.184853553535E+00}
& \num{7.031065924299E+00}
\\
  \runid
&
& \num{5.000000000000E-01}
& \num{1.500000000000E+00}
& \num{6.335631193451E-01}
& \num{7.756824475801E-01}
& \num{2.050882658359E-01}
& \num{2.365005192108E-01}
& \num{3.023152893860E-01}
& \num{1.390993968358E+00}
& \num{2.220956565338E+01}
\\
  \runid
&
& \num{5.000000000000E-01}
& \num{1.000000000000E+00}
& \num{6.579410470744E-01}
& \num{7.936844898501E-01}
& \num{1.801661474559E-01}
& \num{2.235202492212E-01}
& \num{2.268601520672E-01}
& \num{1.802993802047E+00}
& \num{5.780244896485E+01}
\\
  \color{gray}\runid
&
& \color{gray}\num{5.000000000000E-01}
& \color{gray}\num{5.000000000000E-01}
& \color{gray}\num{7.225806451613E-01}
& \color{gray}\num{8.389513108614E-01}
& \color{gray}\num{1.266874350987E-01}
& \color{gray}\num{1.858774662513E-01}
& \color{gray}\num{-1.377100187486E+02}
& \color{gray}\num{3.301794269191E+01}
& \color{gray}\num{1.268261580513E+04}
\\
  \color{gray}\runid
&
& \color{gray}\num{5.000000000000E-01}
& \color{gray}\num{0.000000000000E+00}
& \color{gray}\num{7.193142857143E-01}
& \color{gray}\num{8.367455463972E-01}
& \color{gray}\num{1.357736240914E-01}
& \color{gray}\num{1.830218068536E-01}
& \color{gray}\num{-2.366367496004E+10}
& \color{gray}\num{9.846649297501E+07}
& \color{gray}\num{8.399104116692E+11}
\\
  \runid
&
& \num{1.000000000000E+00}
& \num{4.000000000000E+00}
& \num{5.220404234841E-01}
& \num{6.859744530163E-01}
& \num{3.486500519211E-01}
& \num{2.959501557632E-01}
& \num{7.439742181251E-01}
& \num{1.077101715140E+00}
& \num{1.544667994735E+00}
\\
  \runid
&
& \num{1.000000000000E+00}
& \num{3.000000000000E+00}
& \num{5.536031589339E-01}
& \num{7.126699707714E-01}
& \num{3.149013499481E-01}
& \num{2.720664589823E-01}
& \num{6.215895560324E-01}
& \num{1.095421801373E+00}
& \num{2.223266295366E+00}
\\
  \runid
&
& \num{1.000000000000E+00}
& \num{2.000000000000E+00}
& \num{5.950704225352E-01}
& \num{7.461368653422E-01}
& \num{2.533748701973E-01}
& \num{2.541536863967E-01}
& \num{4.304899436051E-01}
& \num{1.175137410286E+00}
& \num{6.321217197271E+00}
\\
  \runid
&
& \num{1.000000000000E+00}
& \num{1.500000000000E+00}
& \num{6.212929379134E-01}
& \num{7.664166337676E-01}
& \num{2.167705088266E-01}
& \num{2.440290758048E-01}
& \num{3.191943691565E-01}
& \num{1.343546151923E+00}
& \num{1.672189408955E+01}
\\
  \runid
&
& \num{1.000000000000E+00}
& \num{1.000000000000E+00}
& \num{6.504918032787E-01}
& \num{7.882399682161E-01}
& \num{1.876947040498E-01}
& \num{2.274143302181E-01}
& \num{2.553240900307E-01}
& \num{1.661398902047E+00}
& \num{4.184448345338E+01}
\\
  \color{gray}\runid
&
& \color{gray}\num{1.000000000000E+00}
& \color{gray}\num{5.000000000000E-01}
& \color{gray}\num{7.049852037332E-01}
& \color{gray}\num{8.269692923899E-01}
& \color{gray}\num{1.404465212876E-01}
& \color{gray}\num{1.960020768432E-01}
& \color{gray}\num{-9.660727158675E+00}
& \color{gray}\num{1.355719034144E+01}
& \color{gray}\num{3.288990833760E+03}
\\
  \color{gray}\runid
&
& \color{gray}\num{1.000000000000E+00}
& \color{gray}\num{0.000000000000E+00}
& \color{gray}\num{7.205882352941E-01}
& \color{gray}\num{8.376068376068E-01}
& \color{gray}\num{1.298026998962E-01}
& \color{gray}\num{1.858774662513E-01}
& \color{gray}\num{-2.168778286855E+07}
& \color{gray}\num{8.550640115850E+03}
& \color{gray}\num{1.087913114088E+07}
\\
  \color{gray}\runid
&
& \color{gray}\num{1.500000000000E+00}
& \color{gray}\num{0.000000000000E+00}
& \color{gray}\num{7.071167883212E-01}
& \color{gray}\num{8.284339925174E-01}
& \color{gray}\num{1.381100726895E-01}
& \color{gray}\num{1.952232606438E-01}
& \color{gray}\num{-1.024469119288E+02}
& \color{gray}\num{7.604899613979E+01}
& \color{gray}\num{2.561464902217E+04}
\\
  \runid
&
& \bfseries\num{2.000000000000E+00}
& \bfseries\num{0.000000000000E+00}
& \bfseries\num{6.755972315249E-01}
& \bfseries\num{8.063957361759E-01}
& \bfseries\num{1.627725856698E-01}
& \bfseries\num{2.144340602285E-01}
& \bfseries\num{1.500909453354E-01}
& \bfseries\num{3.672517792971E+00}
& \bfseries\num{2.625167665063E+02}
\\
  \runid
&
& \num{3.000000000000E+00}
& \num{0.000000000000E+00}
& \num{5.893608074012E-01}
& \num{7.416324910702E-01}
& \num{2.346832814123E-01}
& \num{2.723260643821E-01}
& \num{2.521236492995E-01}
& \num{1.348755732599E+00}
& \num{1.113948543542E+01}
\\
  \runid
&
& \num{4.000000000000E+00}
& \num{0.000000000000E+00}
& \num{5.604485219164E-01}
& \num{7.183172197544E-01}
& \num{2.733644859813E-01}
& \num{2.863447559709E-01}
& \num{4.561947609615E-01}
& \num{1.162617651199E+00}
& \num{3.644318296816E+00}
\\\midrule
\multicolumn{11}{c}{\tt PROPOSED}
\\\midrule
  \runid
& $H^2$
& \num{4.375000000000E-04} 
& ---
& \num{6.289991796555E-01}
& \num{7.722522976205E-01}
& \num{2.658359293873E-01}
& \num{2.037902388370E-01}
& \num{5.156940681029E-01}
& \num{1.145337079641E+00}
& \num{3.541434521439E+00}
\\
  \runid
&
& \num{7.750000000000E-05} 
& ---
& \num{6.972120164253E-01}
& \num{8.215968419712E-01}
& \num{2.011941848390E-01}
& \num{1.625129802700E-01}
& \num{2.280862700818E-01}
& \num{1.279852308304E+00}
& \num{7.981490327898E+00}
\\
  \runid
&
& \bfseries\num{3.250000000000E-05} 
& ---
& \bfseries\num{7.256381798002E-01}
& \bfseries\num{8.410084898379E-01}
& \bfseries\num{1.695223260644E-01}
& \bfseries\num{1.513499480789E-01}
& \bfseries\num{1.442871176205E-01}
& \bfseries\num{1.381337032546E+00}
& \bfseries\num{1.317457221435E+01}
\\
  \runid
& LS
& \num{2.125000000000E-02} 
& \num{1.000000000000E-01}
& \num{5.698565533504E-01}
& \num{7.259982475904E-01}
& \num{3.211318795431E-01}
& \num{2.471443406023E-01}
& \num{9.477315529108E-01}
& \num{1.055619210531E+00}
& \num{1.161703830432E+00}
\\
  \runid
&
& \num{5.500000000000E-03}
& \num{1.000000000000E-02}
& \num{6.324699163777E-01}
& \num{7.748625687156E-01}
& \num{2.728452751817E-01}
& \num{1.949636552440E-01}
& \num{8.997668935629E-01}
& \num{1.067644629017E+00}
& \num{1.310764235280E+00}
\\
  \runid
&
& \num{4.937500000000E-03}
& \num{1.000000000000E-03}
& \num{6.650667779633E-01}
& \num{7.988469733049E-01}
& \num{2.440290758048E-01}
& \num{1.726375908619E-01}
& \num{6.456085174796E-01}
& \num{1.091783524610E+00}
& \num{1.941117887406E+00}
\\
  \bfseries\runid
&
& \bfseries\num{4.937500000000E-03}
& \bfseries\num{1.000000000000E-04}
& \bfseries\num{7.025398191993E-01}
& \bfseries\num{8.252844500632E-01}
& \bfseries\num{2.061266874351E-01}
& \bfseries\num{1.526479750779E-01}
& \bfseries\num{2.790969771575E-01}
& \bfseries\num{1.156081496345E+00}
& \bfseries\num{4.091922218576E+00}
\\\bottomrule
\end{tabular}
\end{scriptsize}
\end{table}

\setcounter{runidnum}{0}
\begin{sidewaystable}
\caption
{
Quantification of registration quality across multiple datasets. We compare registration quality between the {\tt DEMONS} algorithm and our {\tt PROPOSED} algorithm. We register all remaining {\tt NIREP} datasets ({\tt na02}--{\tt na16}) to the dataset {\tt na01}. We report results with the following settings (from left to right): ($i$) {\tt DEMONS}, $\sigma_u = 0$, $\sigma_d = 1$; ($ii$) {\tt DEMONS}, $\sigma_u = 0$, $\sigma_d = 2$; ($iii$) {\tt DEMONS}, $\sigma_u = 2$, $\sigma_d = 0$; ($iv$) {\tt DEMONS}, $\sigma_u = 4$, $\sigma_d = 0$; ($v$) {\tt PROPOSED}, linear Stokes regularization, $\beta_v=\num{4.937500000000E-03}$, $\beta_w=\num{1E-4}$. These values are based on the experiments in~\tabref{t:performance-evaluation-nirep-parametersweep}; the parameters for {\tt DEMONS} were chosen according to \tabref{t:performance-evaluation-nirep-parametersweep} to deliver results that are consistent with our formulation in terms of the DSC scores or the values for the determinant of the deformation gradient. We report values for the DSC and the $\max$ and $\min$ values for the  determinant of the deformation gradient. We observe that ($i$) when  {\tt DEMONS} outperforms our method in terms of the DSC score, the resulting map gets close to being nondiffeomorphic and ($ii$) when {\tt DEMONS} delivers a map that has a similar quality as the one obtained by our method, its DSC score is lower than ours.
\label{t:performance-evaluation-nirep}
}
\centering
\begin{ssmall}
\begin{tabular}{rrgrrrgrrrgrrrgrrrgrr}
\toprule
& \mcol{4}{{\tt DEMONS} ($\sigma_u=0$; $\sigma_d=1$)}
& \mcol{4}{{\tt DEMONS} ($\sigma_u=0$; $\sigma_d=2$)}
& \mcol{4}{{\tt DEMONS} ($\sigma_u=2$; $\sigma_d=0$)}
& \mcol{4}{{\tt DEMONS} ($\sigma_u=4$; $\sigma_d=0$)}
& \mcol{4}{{\tt PROPOSED}}
\\\midrule
&
& DSC
& $\min$
& $\max$
&
& DSC
& $\min$
& $\max$
&
& DSC
& $\min$
& $\max$
&
& DSC
& $\min$
& $\max$
&
& DSC
& $\min$
& $\max$
\\\midrule
 {\tt na02}
& \runid
& \num{8.015852047556E-01}
& \num{2.174175403328E-01}
& \num{6.063498729499E+01}
& \runid
& \num{7.443433029909E-01}
& \num{4.224464519978E-01}
& \num{7.123601060829E+00}
& \runid
& \num{8.063957361759E-01}
& \num{1.500909453354E-01}
& \num{2.625167665063E+02}
& \runid
& \num{7.183172197544E-01}
& \num{4.561947609615E-01}
& \num{3.644318296816E+00}
& \runid
& \bfseries\num{8.252844500632E-01}
& \num{2.790969771575E-01}
& \num{4.091922218576E+00}
\\
  {\tt na03}
& \runid
& \num{7.854486285449E-01}
& \num{1.675532178158E-01}
& \num{2.997748715939E+01}
& \runid
& \num{7.269763651182E-01}
& \num{4.399891017704E-01}
& \num{3.867686971102E+00}
& \runid
& \num{7.934769947583E-01}
& \num{1.218127802367E-01}
& \num{3.648188405774E+01}
& \runid
& \num{7.045164286544E-01}
& \num{5.061600252726E-01}
& \num{2.522086216325E+00}
& \runid
& \bfseries\num{8.014141167865E-01}
& \num{2.817056008387E-01}
& \num{2.807353508219E+00}
\\
  {\tt na04}
& \runid
& \num{8.112069984632E-01}
& \num{2.592543372322E-01}
& \num{3.550977772006E+01}
& \runid
& \num{7.542482128208E-01}
& \num{4.764574735520E-01}
& \num{4.498951821465E+00}
& \runid
& \bfseries\num{8.172509336225E-01}
& \num{1.889249250079E-01}
& \num{2.076219786143E+02}
& \runid
& \num{7.218465842854E-01}
& \num{5.415278195121E-01}
& \num{4.598451250065E+00}
& \runid
& \num{7.955747955748E-01}
& \num{3.088062317872E-01}
& \num{1.784124725323E+01}
\\
  {\tt na05}
& \runid
& \bfseries\num{8.198726750718E-01}
& \num{1.813195820898E-01}
& \num{3.532881247046E+01}
& \runid
& \num{7.625953244156E-01}
& \num{3.907209215069E-01}
& \num{3.325642513808E+00}
& \runid
& \num{8.183971083136E-01}
& \num{1.275029796413E-01}
& \num{4.099645836556E+01}
& \runid
& \num{7.272500623286E-01}
& \num{3.864260605574E-01}
& \num{2.097359139644E+00}
& \runid
& \num{8.086191508351E-01}
& \num{2.261436786356E-01}
& \num{2.190411734754E+00}
\\
  {\tt na06}
& \runid
& \num{8.172690763052E-01}
& \num{2.227096731886E-01}
& \num{1.669412834054E+02}
& \runid
& \num{7.588807785888E-01}
& \num{4.389820579072E-01}
& \num{4.510445243979E+00}
& \runid
& \bfseries\num{8.202332657201E-01}
& \num{1.772537030657E-01}
& \num{1.834838245205E+02}
& \runid
& \num{7.192961311317E-01}
& \num{3.953687410708E-01}
& \num{4.054614051238E+00}
& \runid
& \num{7.880518862311E-01}
& \num{4.216500459318E-01}
& \num{3.433701477983E+00}
\\
  {\tt na07}
& \runid
& \num{8.318715256331E-01}
& \num{2.370253818108E-01}
& \num{4.262446767902E+01}
& \runid
& \num{7.742176785937E-01}
& \num{4.986397259652E-01}
& \num{4.442286503576E+00}
& \runid
& \bfseries\num{8.371007371007E-01}
& \num{1.162027453706E-01}
& \num{1.920384525398E+01}
& \runid
& \num{7.531215230560E-01}
& \num{5.290781742248E-01}
& \num{3.317859684911E+00}
& \runid
& \num{8.142422706538E-01}
& \num{4.165982937329E-01}
& \num{2.560123277195E+00}
\\
  {\tt na08}
&  \runid
& \num{8.244541484716E-01}
& \num{1.555727243377E-01}
& \num{1.015201662746E+02}
& \runid
& \num{7.635970704766E-01}
& \num{3.742509307936E-01}
& \num{3.523511592092E+00}
& \runid
& \bfseries\num{8.279624460797E-01}
& \num{1.230874400471E-01}
& \num{1.057732761624E+02}
& \runid
& \num{7.305015598752E-01}
& \num{4.341174065432E-01}
& \num{3.337435184733E+00}
& \runid
& \num{8.065017329987E-01}
& \num{3.270254128926E-01}
& \num{2.438659253935E+00}
\\
  {\tt na09}
& \runid
& \bfseries\num{7.990126939351E-01}
& \num{1.358076811859E-01}
& \num{4.726021717033E+01}
& \runid
& \num{7.428105716614E-01}
& \num{3.659671881818E-01}
& \num{7.725577063527E+00}
& \runid
& \num{7.943562610229E-01}
& \num{7.233881781434E-02}
& \num{1.927654172996E+02}
& \runid
& \num{7.183342972817E-01}
& \num{3.621089870222E-01}
& \num{4.374917933613E+00}
& \runid
& \num{7.929400615968E-01}
& \num{2.733781627835E-01}
& \num{5.294247878427E+00}
\\
  {\tt na10}
& \runid
& \num{8.116956141447E-01}
& \num{2.230866959415E-01}
& \num{7.605436521043E+01}
& \runid
& \num{7.546798029557E-01}
& \num{4.634840853298E-01}
& \num{8.398767501250E+00}
& \runid
& \bfseries\num{8.118109238373E-01}
& \num{2.125796796222E-01}
& \num{2.032810552089E+02}
& \runid
& \num{7.268213975425E-01}
& \num{5.301564943008E-01}
& \num{4.089082453133E+00}
& \runid
& \num{7.789060696977E-01}
& \num{3.622765624028E-01}
& \num{3.090187886615E+00}
\\
  {\tt na11}
& \runid
& \num{7.873664435834E-01}
& \num{1.901678208926E-01}
& \num{4.911506035831E+01}
& \runid
& \num{7.271668219944E-01}
& \num{4.225521831015E-01}
& \num{4.612286296842E+00}
& \runid
& \bfseries\num{8.029439696106E-01}
& \num{1.080446667553E-01}
& \num{1.329971249943E+02}
& \runid
& \num{7.032378290240E-01}
& \num{4.291904324841E-01}
& \num{3.354628761757E+00}
& \runid
& \num{7.740058195926E-01}
& \num{2.704628916948E-01}
& \num{1.010060268364E+01}
\\
  {\tt na12}
& \runid
& \num{8.053035589672E-01}
& \num{1.094721888568E-01}
& \num{6.478994110642E+01}
& \runid
& \num{7.502895529303E-01}
& \num{2.979239437738E-01}
& \num{4.084359672418E+00}
& \runid
& \num{7.924924224761E-01}
& \num{6.773363916698E-02}
& \num{8.176397234878E+01}
& \runid
& \num{7.340376665136E-01}
& \num{3.892260862336E-01}
& \num{3.696280033622E+00}
& \runid
& \bfseries\num{8.152576826498E-01}
& \num{4.784719210086E-01}
& \num{3.171984546152E+00}
\\
  {\tt na13}
& \runid
& \bfseries\num{8.254430687060E-01}
& \num{1.623955568344E-01}
& \num{3.611289394948E+01}
& \runid
& \num{7.743352037082E-01}
& \num{3.999776750110E-01}
& \num{3.580397407583E+00}
& \runid
& \num{8.199975789856E-01}
& \num{8.564551320647E-02}
& \num{2.844035710373E+01}
& \runid
& \num{7.499692836958E-01}
& \num{3.992559834995E-01}
& \num{2.818766227567E+00}
& \runid
& \num{8.052930056711E-01}
& \num{2.228257373798E-01}
& \num{2.857687765942E+00}
\\
  {\tt na14}
& \runid
& \bfseries\num{8.068506184586E-01}
& \num{1.415505966648E-01}
& \num{9.625295031322E+01}
& \runid
& \num{7.536892621476E-01}
& \num{3.664275725199E-01}
& \num{6.504831054426E+00}
& \runid
& \num{8.063792085056E-01}
& \num{1.233225203023E-01}
& \num{3.432861699221E+02}
& \runid
& \num{7.328408007626E-01}
& \num{3.858512206079E-01}
& \num{6.434958058506E+00}
& \runid
& \num{7.937355833435E-01}
& \num{2.264179530943E-01}
& \num{9.380807260276E+00}
\\
  {\tt na15}
& \runid
& \bfseries\num{8.148239349914E-01}
& \num{2.548130818281E-01}
& \num{1.698355000634E+02}
& \runid
& \num{7.565313381577E-01}
& \num{4.941050161450E-01}
& \num{9.109601870853E+00}
& \runid
& \num{8.083780077297E-01}
& \num{1.788802307523E-01}
& \num{3.088660718797E+02}
& \runid
& \num{7.333415142962E-01}
& \num{5.050316002037E-01}
& \num{3.346510181087E+00}
& \runid
& \num{7.719015865716E-01}
& \num{4.444753849866E-01}
& \num{1.362253558156E+01}
\\
  {\tt na16}
& \runid
& \bfseries\num{8.093988549618E-01}
& \num{1.663881724866E-01}
& \num{7.377949588274E+01}
& \runid
& \num{7.683467017838E-01}
& \num{4.248637978754E-01}
& \num{5.236314779599E+00}
& \runid
& \num{8.085156993340E-01}
& \num{1.036513074527E-01}
& \num{4.635955407903E+01}
& \runid
& \num{7.448028673835E-01}
& \num{3.966375180016E-01}
& \num{4.386862532188E+00}
& \runid
& \num{8.086000247127E-01}
& \num{3.489095784484E-01}
& \num{5.726002711374E+00}
\\\bottomrule
\end{tabular}
\end{ssmall}
\end{sidewaystable}

\ipoint{Observations} The most important observation is that our framework allows us to generate diffeomorphic maps that are much better behaved with a higher data fidelity, at the cost of a significant increase in computational work as compared to the {\tt DEMONS} algorithm.

The {\tt DEMONS} algorithm is much more efficient than our approach. The time to solution is significantly faster than for our prototype implementation.\footnote{The runtime of our solver depends on the choice of the regularization norm and weight, and the complexity of the registration problem. Considering the test problem in \tabref{t:performance-evaluation-nirep-parametersweep} we obtain the following timings for a fixed regularization parameter on a Linux machine with Intel Xeon X5650 Westmere EP 6-core processors at 2.67GHz with 24GB DDR3-1333 memory (stopping condition: reduction of the gradient by two orders of magnitude): 15 minutes for $\beta_v=\num{1E-2}$, 60 minutes for $\beta_v=\num{1E-3}$, and 300 minutes for $\beta_v=\num{1E-4}$ for an $H^2$-regularization model (compressible diffeomorphism), respectively; 1 minute for $\beta_v = \num{5E-1}$, $\beta_w=\num{1E-4}$, 15 minutes for $\beta_v = \num{5E-2}$, $\beta_w=\num{1E-4}$, and 75 minutes for $\beta_v = \num{5E-3}$, $\beta_w=\num{1E-4}$, for the $H^1$ Stokes regularization scheme, respectively. Applying the {\tt DEMONS} scheme takes only a few seconds. We note that our solver is not finalized; we are currently working on an improved solver in an effort to make it competitive with existing software for diffeomorphic image registration. We will report a detailed performance analysis for this solver elsewhere.} Overall, we trade numerical accuracy and convergence guarantees against computational complexity. The obtained maps are for most of the combinations for $\sigma_u$ and $\sigma_d$ diffeomorphic as judged by the determinant of the deformation gradient. The highest diffeomorphic DSC we could achieve for the {\tt DEMONS} algorithm is \num{7.986184909671E-01} for the \iname{diffusive regularization} (\runref{5} in \tabref{t:performance-evaluation-nirep-parametersweep}) and \num{8.063957361759E-01} for the \iname{fluid-like regularization} (\runref{22} in \tabref{t:performance-evaluation-nirep-parametersweep}). The $\min$/$\operatorname{mean}$/$\max$ values for the determinant of the deformation gradient are \num{2.174175403328E-01}/\num{1.826167477435E+00}/\num{6.063498729499E+01} and \num{1.500909453354E-01}/\num{3.672517792971E+00}/\num{2.625167665063E+02}, respectively. Our formulation allows us to obtain similar (\runref{30} in \tabref{t:performance-evaluation-nirep-parametersweep}) or even better values for the DSC (\runref{26}, \runref{27}, and \runref{31} in \tabref{t:performance-evaluation-nirep-parametersweep}) with much more well-behaved deformation maps (as judged by the determinant of the deformation gradient).

For \runref{30} we obtain almost the same DSC with $\det(\mat{F}_1) \in [\num{6.456085174796E-01},\num{1.941117887406E+00}]$. Further, we can observe that across almost all runs the mean values for $\det(\mat{F}_1)$ are much closer to one for the linear Stokes regularization case. The $H^2$-regularization also results in better behaved deformation maps. However, the variations increase significantly as we turn to a higher data fidelity. Our approach results in maps for which the maximal value of the determinant of the deformation gradient is (for the most part) much better behaved than for the {\tt DEMONS} algorithm. This is important since the Jacobian of the inverse deformation map will have very small values if the maximum for the reported values is large. For instance, for the best runs for {\tt DEMONS}  $\max(\det(\mat{F}_1))$ is equal to \num{6.063498729499E+01} (pure \iname{diffusive regularization}; \runref{5} in \tabref{t:performance-evaluation-nirep-parametersweep}) and \num{2.625167665063E+02} (pure \iname{fluid-like regularization}; \runref{22} in \tabref{t:performance-evaluation-nirep-parametersweep}) as compared to \num{4.091922218576E+00} for the linear Stokes case (\runref{31} in \tabref{t:performance-evaluation-nirep-parametersweep}; which in addition to that has a better DSC). If we compare the results with similar values for the deformation gradient (e.g., \runref{3} or \runref{24} in \tabref{t:performance-evaluation-nirep-parametersweep}) we cannot achieve the same DSC scores as the {\tt PROPOSED} algorithm delivers; we have to operate the {\tt DEMONS} algorithm at regimes with large variations in the determinant of the deformation gradient to obtain DSC scores that are equivalent to those achieved with our algorithm.

If we consider the results for the remaining datasets in \tabref{t:performance-evaluation-nirep} we observe a similar behavior (for fixed parameters for both algorithms). We obtain slightly better DSC scores for the {\tt DEMONS} algorithm (see also \figref{f:results-nirep-quantitative-evaluation-boxplots}) at the cost of a larger variations in the determinant of the deformation gradient. If we increase the regularization we can reproduce similar values for the determinant of the deformation gradient but are not able to achieve the same data fidelity as judged by the DSC scores. \figref{f:results-nirep-quantitative-evaluation-boxplots} shows that these differences are on average consistent across all datasets.

\ipoint{Conclusions} We have conducted a preliminary two-dimensional study of registration quality based on the {\tt NIREP} data. All approaches deliver diffeomorphic maps with a good data fidelity. The {\tt DEMONS} algorithm arrives at a solution significantly faster than our current prototype implementation. We note that our method has not been optimized for speed yet; there exist several ways to accelerate our algorithm, which we are currently investigating. We will report these improvements and the extension of our solver to three-dimensional problems elsewhere. This preliminary study suggests that our algorithm provides much more well-behaved mappings without compromising data fidelity as compared to the {\tt DEMONS} algorithm. Overall, we trade numerical accuracy and convergence guarantees against computational efficiency (i.e., an increase in the time to solution). We consider these differences in registration quality a preliminary result; clearly, we have to validate the performance of our solver on three-dimensional data, extend the comparison to other algorithms, and reduce the time to solution to truly demonstrate that our formulation has the potential to impact the applied sciences.

\subsection{Nonlinear Stokes regularization (shear control)}
\label{s:res-shear-control}

\ipoint{Purpose} We study the effect of controlling the shear in the deformation field in the presence of an expected `'discontinuous`' motion field. We compare results for the nonlinear Stokes regularization model to plain, quadratic smoothness regularization, and a linear Stokes regularization model (incompressible diffeomorphism).

\begin{figure}
\centering
\includegraphics[width=0.99\textwidth]
{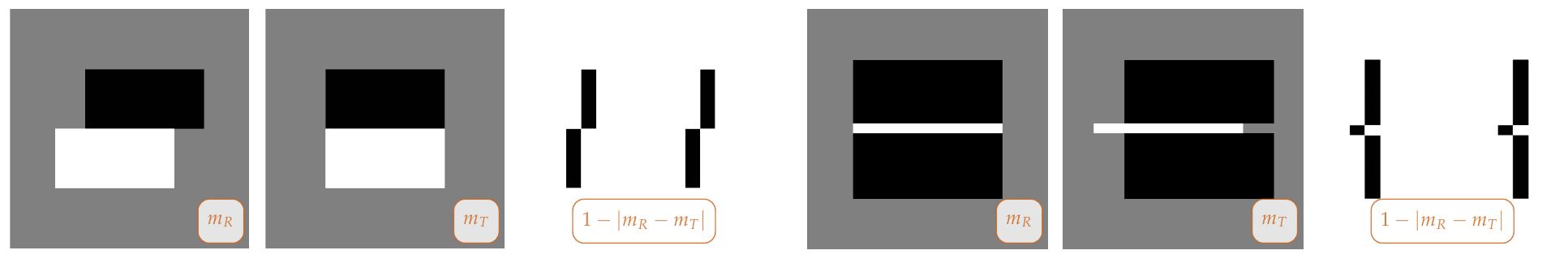}
\caption{Registration problems with an expected `'discontinuous`' motion field (sliding interfaces; left: \iname{sliding rectangles}; right: \iname{sliding vent}). We display (from left to right) the reference image $m_R$ (fixed image), the template image $m_T$ (image to be registered), and a map of the residual differences between $m_R$ and $m_T$ before registration (for each set of images as indicated by the inset).}
\label{f:regprob-image-collection-sliding}
\end{figure}

\ipoint{Setup} We consider two synthetic problems for which we expect the deformation to contain large shear (see \figref{f:regprob-image-collection-sliding}). The images have a grid size of $512\times512$.\footnote{We use more grid points to be able to resolve the velocity field and avoid aliasing.} We compare plain smoothness regularization based on an $H^2$-seminorm ($\gamma=0$) to models of incompressible flow ($H^1$-regularization; $\gamma=1$). We study the qualitative behavior of the deformation map with respect to changes in the flow law exponent $\nu$ for empirically chosen values for $\beta_v\in\{\num{1E-2}, \num{1E-3}\}$. In particular, we study shearthickening ($\nu=1/2$) and shearthinning ($\nu\in\{3,5\}$). We consider the full set of termination criteria used in~\cite{Mang:2015a} for this set of experiments with a tolerance of $\num{1E-3}$. No grid, scale or parameter continuation is performed.

\begin{figure}
\centering
\includegraphics[width=0.99\textwidth]
{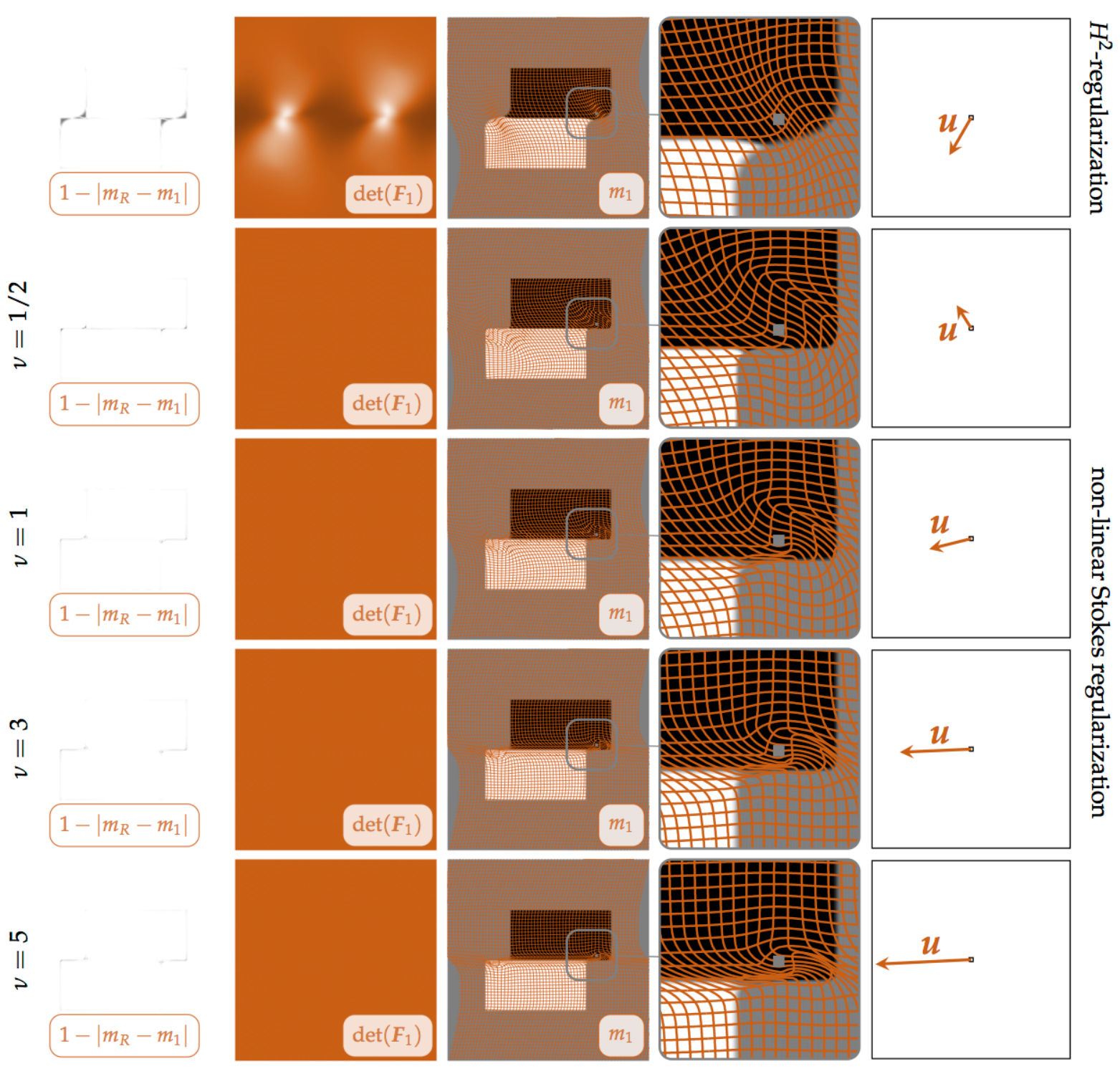}
\caption
{
Exemplary registration results for the \iname{sliding rectangles} (see \figref{f:regprob-image-collection-sliding}). We study the effect of shear control (nonlinear Stokes regularization). We compare plain $H^2$-regularization (top row; $\gamma=0$; $\beta_v=\num{1E-2}$) to a linear Stokes regularization model (third row from the top; $H^1$-regularization; $\gamma=1$; $\beta_v=\num{1E-3}$) and a nonlinear Stokes regularization model (second row from the top: $\nu=1/2$ (shear thickening; $\beta_v=\num{1E-2}$); first and second row from the bottom: $\nu\in\{3,5\}$ (shear thinning; $\beta_v=\num{1E-3}$)). We show (from left to right) ($i$) a map of the residual differences between the reference image $m_R$ and the deformed template $m_1$, ($ii$) a map of the determinant of the deformation gradient, ($iii$) a deformed grid overlaid onto the deformed template image $m_1$ (to illustrate the deformation map $\vect{y}$), ($iv$) a close up of the latter for a particular area of interest, and ($v$) a single displacement vector at $\vect{x} = (\num{4.663302e+00},\num{3.252039e+0})$ (the location is indicated as a gray rectangle in the visualization of the deformed grid; the size of the box is $25\times25$ grid points).
}
\label{f:res-sliding-rectangle}
\end{figure}

\begin{figure}
\centering
\includegraphics[width=0.99\textwidth]
{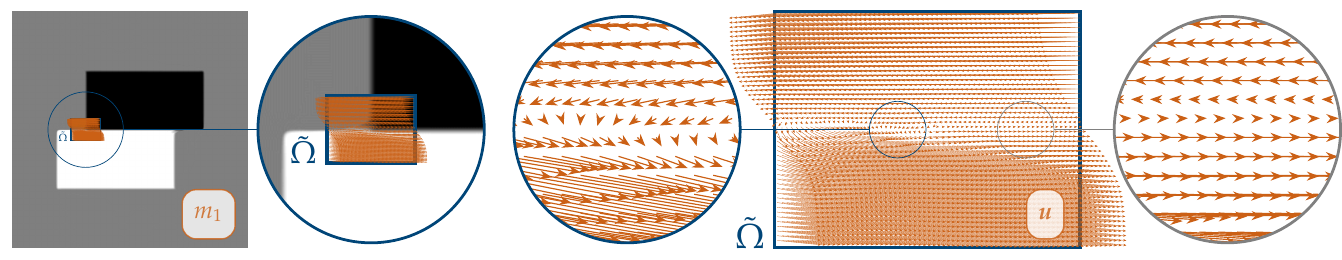}
\caption
{
Exemplary registration results for the \iname{sliding rectangles} (see \figref{f:regprob-image-collection-sliding}). We display the displacement field for the nonlinear Stokes regularization for $\nu=5$ (bottom row in \figref{f:res-sliding-rectangle}). We only show an detail of the displacement field in full resolution. On the left we illustrate the region of interest $\tilde{\Omega}$. A closeup of this region of interest is provided on the right.
}
\label{f:res-sliding-rectangle-displacement}
\end{figure}

\begin{figure}
\centering
\includegraphics[width=0.9\textwidth]
{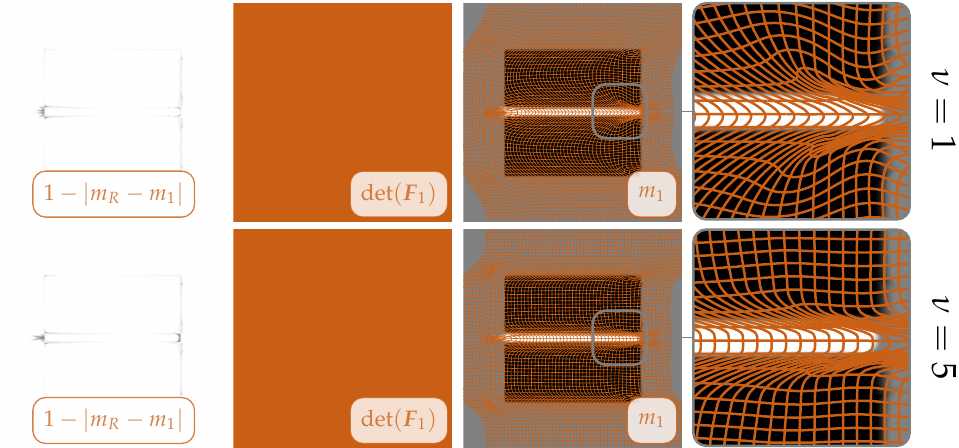}
\caption
{
Exemplary registration results for the \iname{sliding vent} (see \figref{f:regprob-image-collection-sliding}). We study the effect of shear control for a highly nonlinear registration problem. We show results for a linear (top row; $\gamma = 1$; $\nu=1$) and a nonlinear (bottom row; $\gamma=1$; $\nu=5$) Stokes regularization model. We show (from left to right) ($i$) a map of the residual differences between the reference image $m_R$ and the deformed template $m_1$, ($ii$) a map of the determinant of the deformation gradient, ($iii$) a deformed grid overlaid onto the deformed template image $m_1$ (to illustrate the deformation map $\vect{y}$) and ($iv$) a close up of the latter for a particular area of interest.
}
\label{f:res-sliding-vent}
\end{figure}

\ipoint{Results} We report exemplary results for the \iname{sliding rectangles} and the \iname{sliding vent} in \figref{f:res-sliding-rectangle}, \figref{f:res-sliding-rectangle-displacement}, and \figref{f:res-sliding-vent}, respectively. We enforce incompressibility up to numerical accuracy.

\ipoint{Observations} The most important observation is that the nonlinear Stokes regularization provides an adaptive control of the shear of the deformation field at the sliding interface. Setting $\nu$ to a value in $(0,1)$ increases the resistance to shear (shear thickening fluid). On the contrary, if we choose $\nu>1$ we promote shear. The larger $\nu$ the sharper the transition at the interface and the more localized the deformation. This confirms the theoretical statement that the model tends to a total variation regularization for $\nu\rightarrow\infty$. However, we can already recover sharp interfaces for small $\nu$ (e.g., for $\nu = 5$; see \figref{f:res-sliding-rectangle} and \figref{f:res-sliding-vent} bottom and, in particular, \figref{f:res-sliding-rectangle-displacement}). The residual differences between the registered images are insignificant for varying parameters $\nu$. The computed mappings are very different; points close to the sliding interface map to completely different positions. We can model highly nonlinear deformations with a precise control on the determinant of the deformation gradient. Likewise to the linear case we can also extend this formulation by introducing a mass source $w$ (see \secref{s:derivation-non-linear-stokes-model} in the appendix for details) rendering the flow near-incompressible (results not included in this study). This makes this approach applicable across a wider range of registration scenarios.

\ipoint{Conclusion} The nonlinear Stokes regularization model allows us to promote or penalize large shear in the deformation field as required. As a consequence, we have---in contrast to traditional vectorial total variation---complete control on the smoothness properties of the deformation map and the determinant of the deformation gradient. Further, we---like in total variation regularization---do not have to identify the interfaces where the sliding motion is expected to occur (i.e.\ we do not require a presegmentation of the data). We note that promoting shear is only an approximation of true sliding motion, i.e., our formulation does not allow for computing `'discontinuous`' motion fields. In future studies we will compare our formulation against total variation regularization to better quantify the capabilities.

\section{Conclusions}
\label{s:conclusions}

We have introduced novel constrained regularization schemes for large deformation diffeomorphic image registration that feature a local control of the divergence of the velocity and thus of the determinant of the deformation gradient (in a problem-dependent way). Our formulation is founded on well-established computational models in fluid mechanics (Stokes flow). All results reported in this study are limited to the two-dimensional case. Nothing in our formulation is specific to the two-dimensional case; an extension to three dimensions is ongoing work in our group.

We invert for a stationary velocity field. We achieve a similar or even better inversion quality (as judged by values for the residual differences and overlap between anatomical regions) as compared to available diffeomorphic registration models, while maintaining a better control over the deformation regularity. Furthermore, several applications do require incompressible or near-incompressible deformations, for example, in medical imaging. Our framework provides such a technology.

It is unclear how to theoretically determine the behavior of the proposed methods. For this reason we conducted experiments to probe their behavior. The basic conclusions from our experiments are the following:
\begin{itemize}
\item Using an $H^1$-seminorm as a regularization model without control of the determinant of the deformation gradient is not robust. Either it produces uninformative maps (large regularization) or highly (perhaps unacceptably so) deformed maps.
\item The $H^2$-seminorm without control of the determinant of the deformation gradient behaves well but its cost with decreasing regularization increases to regimes that make it not practical.
\item Our proposed $H^1$-seminorm regularization plus control of the deformation gradient performs well. It delivers small mismatch values (comparable to the $H^2$ case without controlling $\det(\igrad\vect{y})$) and smooth $\det(\igrad\vect{y})$ much faster than the $H^2$ scheme (e.g., \tabref{t:quantitative-results-quadratic-regularization}, \runref{4} versus \runref{16} and \runref{20} versus \runref{23}). It also delivers a good agreement between anatomical structures with an excellent control on $\det(\igrad\vect{y})$ (e.g., \tabref{t:performance-evaluation-nirep-parametersweep} \runref{29} and \runref{30}). However, this scheme results in one additional regularization parameter $\beta_w$ and its calibration is expensive. This cost can be amortized in studies involving multiple images.
\end{itemize}

So, the new formulation seems preferable in terms of robustness and speed. Further studies in three dimensions and on a larger set of images are necessary to confirm these results.

We have also introduced a regularization model that allows us to control the shear in the deformation map. We can either promote (shear thinning) or penalize (shear thickening) shear in an attempt to \emph{approximate} discontinuous motion fields or generate more well-behaved deformation maps. Our framework \bipa\item is applicable to smooth and non-smooth registration problems, \item allows us to control the amount of shear, \item does not require a presegmentation of the data, and \item features a control of the determinant of the deformation gradient in a problem dependent way\eipa. We demonstrated that this regularization results in dramatically different deformation maps (see \figref{f:regprob-image-collection-sliding}).

Our solver is not finalized; its extension to three dimensions requires more work. The next steps will be the design of a more effective preconditioner, the implementation of a more effective scheme to solve the transport equations, and the application to problems that have time sequences of images. For such cases, a time-dependent velocity field will be necessary.

\subsection*{Acknowledgments} We thank Georg Stadler and Johann Rudi for helpful discussions and suggestions.

\begin{appendix}

\section{Variable Elimination: Derivations}
\label{s:constraint-elimination-derivation}

Here, we provide the derivations of the variable elimination. We start with the linear Stokes regularization model (i.e., we consider the regularization operator $\D{A}=-\ilap$ and $\gamma=1$ in \eqref{e:first-order-opt}).

\subsection{Linear Stokes Regularization}
\label{s:derivation-linear-stokes-model}

The variable elimination for the linear, incompressible flow model (i.e.,enforcing $\idiv\vect{v} = 0$ up to numerical accuracy) can be found in~\cite{Mang:2015a}. Here, we extend on the formulation in~\cite{Mang:2015a} by introducing a mass source $w$ into the divergence constraint. We eliminate $p$ and $w$ from the optimality system~\eqref{e:first-order-opt}. This will result in an optimality system that allows us to only iterate on the reduced space of the control variable $\vect{v}$.

Applying the divergence to~\eqref{e:control-eq-v} yields $\idiv(-\beta_v\ilap\vect{v}) + \ilap p + \idiv\vect{b} = 0$. From the optimality condition $\idiv \vect{v} = w$ and the equivalence
\begin{equation}
\label{e:laplacian-div-curl}
\ilap \vect{v} = \igrad(\idiv \vect{v}) - \icurl(\icurl \vect{v})
\end{equation}

\noindent it follows that $-\beta_v\ilap w + \ilap p = -\idiv \vect{b}$. From~\eqref{e:control-eq-w} it follows that
\begin{equation}
\label{e:mass-source}
w = -(\beta_w(-\ilap  + \id))^{-1}p.
\end{equation}

\noindent Inserting this expression yields $\ilap(\beta_v(\beta_w(-\ilap+\id))^{-1} + \id) p = -\idiv\vect{b}$ and therefore
\begin{equation}
\label{e:pressure-elimination-relaxed-ic}
p =
-
(
\beta_v(\beta_w(-\ilap+\id))^{-1} + \id
)^{-1}
\ilap^{-1}\idiv\vect{b}.
\end{equation}

\noindent Inserting this expression into~\eqref{e:control-eq-v} yields the control equation
\begin{equation}
\label{e:derive-control-eq-elim}
-\beta_v\ilap\vect{v}
\underbrace
{
-
\igrad
(
\beta_v(\beta_w(-\ilap+\id))^{-1} + \id
)^{-1}
\ilap^{-1}
\idiv\vect{b}
+ \vect{b}
}_{\eqdef\D{K}[\vect{b}]}
\end{equation}

\begin{figure}
\centering
\includegraphics[width=0.65\textwidth]
{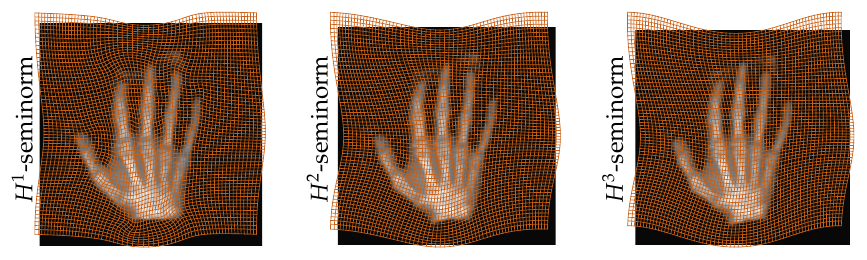}
\caption
{
Exemplary results for a fully incompressible Stokes model (i.e., $w=0$) for different regularization norms. We provide (from left to right) results for an $H^1$-, an $H^2$-, and an $H^3$-seminorm for $\beta_v=\num{1E-2}$.
}
\label{f:results-hands-stokes-different-regularization-norms}
\end{figure}

\noindent Note, that~\eqref{e:derive-control-eq-elim} is independent of the variables $w$ and $p$. In addition, we have eliminated~\eqref{e:ic-constraint} and~\eqref{e:control-eq-w} from~\eqref{e:first-order-opt}. We arrive at the optimality system~\eqref{e:first-order-opt-elim}. Computing second variations of the weak form of~\eqref{e:first-order-opt-elim} yields~\eqref{e:second-order-elim} with the operator $\D{K}$ as defined in~\eqref{e:derive-control-eq-elim}. If we consider an incompressible diffeomorphism (i.e., set $w=0$) the control equation~\eqref{e:derive-control-eq-elim} simplifies to $-\beta_v\ilap\vect{v} -\igrad\ilap^{-1}\idiv \vect{b} + \vect{b} = 0$~\cite{Mang:2015a}. Note that it is possible to replace the Laplacian operator with a biharmonic operator (i.e., consider an $H^2$- instead of an $H^1$-seminorm); the same arguments used above still hold. We can even use an $H^3$-seminorm (the reduced gradient becomes a triharmonic equation) if theoretical considerations are of concern (see~\cite{Chen:2011a} or \secref{s:theoretical-considerations} for a brief discussion; we provide exemplary results in \figref{f:results-hands-stokes-different-regularization-norms}). We limit ourselves in this work to an $H^1$-seminorm, which results in a linear Stokes regularization model.

\subsection{Nonlinear Stokes Regularization}
\label{s:derivation-non-linear-stokes-model}

We discuss the variable elimination techniques for the nonlinear Stokes regularization next. We start with a model of incompressible flow (i.e., we assume that $\idiv\vect{v} = 0$). The same arguments that have been used in the former section apply. However, since the viscosity is no longer a constant but a function of $\vect{v}$, we have to decompose $\eta$ into
\begin{equation*}
\eta[\vect{v}] = \hat{\eta}[\vect{v}] + \bar{\eta}[\vect{v}]
= \hat{\eta}[\vect{v}] + \frac{1}{\#\Omega}\iom{\eta[\vect{v}]}
\end{equation*}

\noindent to be able to eliminate $p$. If we insert this decomposition into the control equation for $\vect{v}$ we obtain
\begin{equation}
\label{e:control-eq-v-ns-decomp}
-2\beta_v\bar{\eta}[\vect{v}]\idiv\D{E}[\vect{v}]
-\idiv2\beta_v\hat{\eta}[\vect{v}]\D{E}[\vect{v}]
+\igrad p + \vect{b}
= 0.
\end{equation}

\noindent The divergence of the strain rate tensor $\D{E}[\vect{v}]$ is identical to $\half{1}\ilap \vect{v}$ under the incompressibility assumption $\idiv \vect{v} = 0$. Accordingly, we have
\[
-\beta_v\bar{\eta}[\vect{v}] \ilap \vect{v}
-\idiv2\beta_v\hat{\eta}[\vect{v}]\D{E}[\vect{v}]
+\igrad p + \vect{b}
= 0.
\]

\noindent By taking advantage of~\eqref{e:laplacian-div-curl} and applying the divergence we obtain
\[
-\idiv\idiv2\beta_v\hat{\eta}[\vect{v}]\D{E}[\vect{v}]
+ \ilap p + \idiv\vect{b} = 0
\]

\noindent and therefore
$
p =
\ilap^{-1}\idiv
(
\idiv2\beta_v\hat{\eta}[\vect{v}]\D{E}[\vect{v}] - \vect{b}
)
$. Inserting this expression into the control equation for $\vect{v}$ results in
\begin{equation}
\label{e:control-eq-ns}
\vect{\tilde{g}}
\defeq
-\idiv2\beta_v\eta[\vect{v}]\D{E}[\vect{v}]
+ \D{K}[\vect{b},\vect{v}]=0,
\end{equation}

\noindent where
$
\D{K}[\vect{b},\vect{v}]
=
\igrad
\ilap^{-1}
\idiv
(
\idiv2\beta_v\hat{\eta}[\vect{v}]\D{E}[\vect{v}]
-
\vect{b}
)
+
\vect{b}
$. Computing second variations yields the incremental control equation
\begin{equation}
\label{e:inc-control-eq-elim-nlstokes}
\beta_v\D{B}(\vect{v})[\vect{\tilde{v}}]
+\D{L}(\vect{v})[\vect{\tilde{b}}, \vect{\tilde{v}}]
=-\vect{\tilde{g}}.
\end{equation}

\noindent The operator $\D{B}$ is the second variation of~\eqref{e:shear-control-reg-v} given in~\eqref{e:2nd-var-reg-v-nls} and
\[
\D{L}(\vect{v})[\vect{\tilde{b}},\vect{\tilde{v}}]
=
\igrad\ilap^{-1}\idiv
\left(
\idiv
2\beta_v
\left(
\hat{\eta}[\vect{v}]
+
\eta[\vect{v}]
\D{Q}[\vect{v}]
\right)
\D{E}[\vect{\tilde{v}}]
-
\vect{b}
\right)
+
\vect{b}.
\]

Next, we consider a nonzero mass source $w$ (i.e., we relax the incompressibility constraint to $\idiv v = w$). In this case, the divergence of the strain rate tensor $\D{E}[\vect{v}] = \half{1}((\igrad\vect{v}) + (\igrad\vect{v})^\T)$ is no longer proportional to $\ilap\vect{v}$. Instead, we have
\[
-\beta_v\bar{\eta}[\vect{v}](\ilap\vect{v} + \igrad w)
-\idiv2\beta_v\hat{\eta}[\vect{v}]\D{E}[\vect{v}]
+\igrad p + \vect{b}
=0.
\]

\noindent Applying the divergence operator yields
\[
-2\beta_v\bar{\eta}[\vect{v}]\ilap w
-\idiv\idiv2\beta_v\hat{\eta}[\vect{v}]\D{E}[\vect{v}]
+\ilap p + \idiv\vect{b}
= 0.
\]

\noindent Using~\eqref{e:mass-source} we can eliminate $w$. Thus,
\[
\ilap
(
\underbrace
{
2\beta_v\bar{\eta}[\vect{v}](\beta_w(-\ilap+\id))^{-1}+\id
}_{\eqdef\D{M}}
)
p
-\idiv\idiv2\beta_v\hat{\eta}[\vect{v}]\D{E}[\vect{v}]
= -\idiv\vect{b}
\]

\noindent and therefore
$
p
=
\D{M}^{-1}
\ilap^{-1}
\idiv
\left(
\idiv2\beta_v\hat{\eta}[\vect{v}]\D{E}[\vect{v}]
-\vect{b}
\right)
$. We have again found an expression for $p$ that not only eliminates $p$ but also the equation for the constraint on $\idiv\vect{v}$, the control variable $w$, and the associated control equation for $w$. We obtain first order optimality conditions that are very similar to the incompressible case. The only difference is the form of the operator $\D{K}$ in~\eqref{e:control-eq-ns}. In particular,
\[
\D{K}[\vect{b},\vect{v}]
=
\igrad
\D{M}^{-1}\ilap^{-1}
\idiv
\left(
\idiv2\beta_v\hat{\eta}[\vect{v}]\D{E}[\vect{v}]
-
\vect{b}
\right)
+
\vect{b}.
\]

\noindent It immediately follows that we arrive at~\eqref{e:project-inc-body-force-nls} for the second variation.

\section{Theoretical Considerations}
\label{s:theoretical-considerations}

An adequate choice for the regularization norm ensures the existence and uniqueness of a minimizer for $\D{J}$ in~\eqref{e:objective}. A second requirement for an admissible solution of~\eqref{e:objective} is that the map $\vect{y}$ generated by $\vect{v}\in\fs{V}$ is a diffeomorphism. In~\cite{Dupuis:1998a,Trouve:1995a} it is shown in the context of the \iname{large deformation diffeomorphic metric mapping} formulation (which is related to our formulation; see~\cite{Hart:2009a,Mang:2015a}) that $\vect{y}$ is diffeomorphic if we enforce a sufficient amount of smoothness on the elements of the space $\fs{V}$. The order of the associated Sobolev norm depends on the dimensionality $d$ of the ambient space; for $d=3$ the norm $\|(-\ilap+\id)^\gamma\vect{v}\|^2_{L^2(\Omega)^d}$, where $\gamma > 1.5$, is adequate~\cite{Beg:2005a}.

Our formulation operates in an incompressible or near-incompressible regime, i.e., $\vect{v}\in\{\vect{f}\in H^1(\Omega)^d : \idiv \vect{f} = w, w \in H^1(\Omega)\}$. An excellent reference for the analysis of such flows is~\cite{Crippa:2007a}. An existence proof for a minimizer of $\fs{J}$ in~\eqref{e:opt-prob-v} for the nonstationary, incompressible case (i.e., $w=0$) in two dimensions, for which the images are modeled as functions of bounded variation and the regularization model for $\vect{v}$ is an $H^3$-seminorm, can be found in~\cite{Chen:2011a}. This proof is based on the premise that $\vect{v}$ propagates $\fs{BV}$ regularity in space (i.e., $m_T\in\fs{BV}(\Omega)$ and $m(\cdot,t)$ is in $\fs{BV}(\Omega)$ for all $t\in[0,1]$); Lipschitz regularity of $\vect{v}$ in space not only implies that we propagate the regularity of $m_T$ to $t\in[0,1]$ but also that we obtain a unique map $\vect{y}$ from $\vect{v}$~\cite{Chen:2011a,Crippa:2007a}. Under these assumptions it is shown in~\cite{Chen:2011a} that we have to equip $\fs{J}$ with an $H^3$-seminorm to guarantee existence of a minimizer; this will result in a triharmonic control equation. The results reported in~\cite{Chen:2011a} were---due to numerical considerations---obtained using an $H^1$-seminorm assuming that $H^1$-regularity still yields smooth enough results. A critical result to support this claim can be found in~\cite{DiPerna:1989a} (see also \cite{Crippa:2007a}); the smoothness requirements for $\vect{v}$ are relaxed to a local Sobolev regularity; this relaxation will not preserve $\fs{BV}$ regularity of the images~\cite{Colombini:2004a}. However, existence and uniqueness of the flow $\vect{y}$ can be guaranteed if $\vect{v}$ has local Sobolev regularity in space and $\idiv \vect{v}$ is a bounded measurable function~\cite{DiPerna:1989a}; our incompressible formulation fulfills these requirements. According to~\cite{DiPerna:1989a}, $H^1$ regularity of $\vect{v}$ will transport $L^2$ images to $L^2$ images (see also~\cite{Chen:2011a,Chen:2011b}). Relaxing $\fs{BV}$-regularity to $L^2$-integrability seems reasonable, especially since our scheme currently cannot handle $\fs{BV}$ images; we have to use smooth representations of $m_T$ and $m_R$ to ensure numerical stability. An existence proof for a minimizer of an incompressible $H^1$-flow that follows these arguments can be found in~\cite[pages~58ff]{Chen:2011b}. A rigorous proof for our formulation remains open. We note that we can switch to an $H^2$- or an $H^3$-seminorm if our formulation does not meet the theoretical requirements; all the derivations and algorithmic features presented here will still apply (see \secref{s:constraint-elimination-derivation}; an exemplary result can be found in \figref{f:results-hands-stokes-different-regularization-norms}). Likewise, adding a parabolic regularization via a diffusion operator to the hyperbolic transport equation can be another strategy to ensure existence of a minimizer, even for an $L^2$-integrable $\vect{v}$~\cite{Barbu:2016a}.

As we mentioned in~\secref{s:problem-formulation}, our computational studies as well as the results reported in~\cite{Chen:2011a} suggest that an $H^1$-seminorm together with a control on $\idiv\vect{v}$ seems to provide sufficient smoothness to converge to a locally optimal, diffeomorphic solution.

\section{Connection to Total Variation Regularization}
\label{s:total-variation}

Here we briefly comment on the connection between our nonlinear Stokes regularization model and a total variation regularization model. The total variation regularization model is given by
\begin{equation}
\label{e:tv}
\F{S}[\vect{v}] = \int_{\Omega}(\igrad\vect{v}:\igrad\vect{v})^{1/2}\d{\vect{x}}
\end{equation}

\noindent Notice that the exponent in~\eqref{e:shear-control-reg-v} will tend to $1/2$ if we let $\nu$ in~\eqref{e:shear-control-reg-v} tend to $\infty$. Thus, if we replace the strain rate tensor $\D{E} = \half{1}((\igrad\vect{v}) + (\igrad\vect{v})^\T)$ with $\igrad\vect{v}$, \eqref{e:shear-control-reg-v} and \eqref{e:tv} are equivalent as $\nu\rightarrow\infty$. This equivalence is also reflected by the variations of both models. We obtain
$
\D{A}[\vect{v}] = -\idiv (\igrad\vect{v}:\igrad\vect{v})^{-1/2} \igrad \vect{v}
$
for the first variation of \eqref{e:tv} and
\[
\D{B}(\vect{v})[\vect{\tilde{v}}]
= -\idiv(\igrad\vect{v}:\igrad\vect{v})^{-1/2}
\left(\D{I}+\frac{\igrad\vect{v}\otimes\igrad\vect{v}}{\igrad\vect{v}:\igrad\vect{v}}\right)
\igrad\vect{\tilde{v}}
\]

\noindent for the second variation of \eqref{e:tv} with respect to $\vect{v}$, respectively. These expressions are very similar to the operators in \eqref{e:1st-var-reg-v-nls} and~\eqref{e:2nd-var-reg-v-nls}. As such, the derivations we presented in this work also hold if we replace \eqref{e:shear-control-reg-v} with \eqref{e:tv}.

\section{Illustration of Deformation Map}
\label{s:visualization}

We report images of the deformation pattern (deformed grid) and maps of the determinant of the deformation gradient in order to illustrate local properties of the deformation map. Here, we provide information on how these were generated and on how to interpret them.

\subsection{Deformation Map}
\label{s:deformation-map}

\begin{figure}
\centering
\includegraphics[width=0.9\textwidth]
{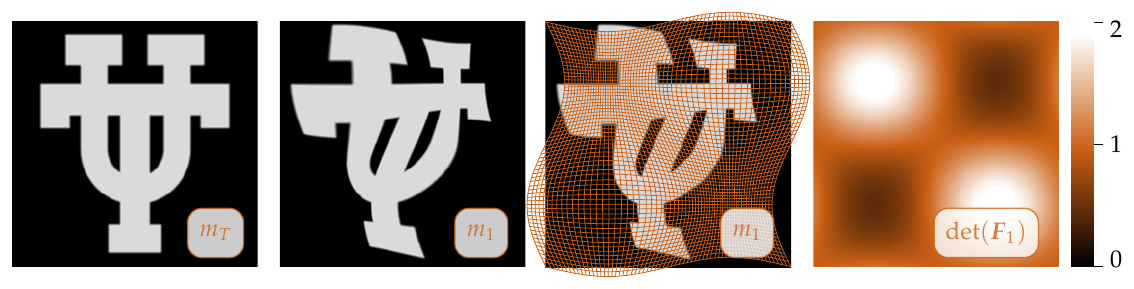}
\caption
{
Illustration of the visualization of the computed deformation map $\vect{y}$. From left to right: template image $m_T$, deformed template image $m_1$ (deformed configuration), deformed template image $m_1$ with an illustration of the deformed grid as an overlay and a map of the determinant of the deformation gradient $\mat{F}_1$ (as identified in the inset). The color map for $\det(\mat{F}_1)$ is displayed on the right. Notice that the map $\vect{y}$ is defined in an Eulerian frame of reference (Eulerian description of motion), i.e., the map $\vect{y}$ models where points originate from. Accordingly, we have $\det(\mat{F}_1)\equiv\det(\igrad\vect{y}^{-1}) = \det(\igrad\vect{y})^{-1}$.
}
\label{f:defgrid-jacobian-illustrated}
\end{figure}

We illustrate regularity and local properties of the deformation map $\vect{y}$ on the basis of deformed grids. We define $\vect{y}$ as a perturbation from identity, i.e.\ $\vect{y} \defeq \vect{x} - \vect{u}_1$, where $\vect{u}_1: \bar{\Omega} \rightarrow \ns{R}^{d}$, $\vect{u}_1 \defeq \vect{u}(\cdot, t=1)$, $\vect{u} :\bar{\Omega} \times [0,1] \rightarrow\ns{R}^d$, is some displacement field at final time $t=1$. The latter can be computed from the velocity field $\vect{v}$ by solving
\begin{equation}\label{e:displacement}
\p_t \vect{u} + (\igrad \vect{u}) \vect{v} = \vect{v}
\;\;{\rm in}\;\; \Omega \times (0,1],
\qquad \vect{u} = 0 \;\;{\rm in}\;\;\Omega \times \{0\},
\end{equation}

\noindent with periodic boundary conditions on $\p\Omega$. Note that $\vect{y}$ is defined in an Eulerian frame of reference (i.e., the deformed illustrates not where the points move to but where they originate from). An exemplary visualization of a synthetic deformation is shown in \figref{f:defgrid-jacobian-illustrated}.

\subsection{Deformation Gradient}
\label{s:deformation-gradient}

We report maps and values for the determinant of the deformation gradient to qualitatively and quantitatively assess regularity of a mapping $\vect{y}$. In the framework of continuum mechanics the deformation tensor field $\mat{F} : \bar{\Omega} \times [0,1] \rightarrow \ns{R}^{d\times d}$ can be computed from $\vect{v}$ by solving
\begin{equation}\label{e:def-grad}
\partial_t \mat{F} + (\vect{v}\cdot \igrad) \mat{F}
= (\igrad \vect{v}) \mat{F} \;\; {\rm in} \;\; \Omega \times (0,1],
\qquad \mat{F} = \mat{I} \;\; {\rm in} \quad \Omega \times \{0\},
\end{equation}

\noindent with periodic boundary conditions on $\p\Omega$. Here, $\mat{I} = \operatorname{diag}(1,\ldots,1)\in\ns{R}^{d\times d}$ and $\det(\mat{F}_1)$ is identical to $\det(\igrad \vect{y})^{-1}$, where $\mat{F}_1 \defeq \mat{F}(\cdot, t=1)$, $\mat{F}_1 : \bar{\Omega} \rightarrow \ns{R}^{d \times d}$, and $\vect{y}$ is the Eulerian deformation map.

We limit the color map for the display of $\det(\mat{F}_1)$ to $[0,2]$. In particular, the color map ranges from black (compression: $\det(\mat{F}_1) \in (0,1)$; black corresponds to values of 0 or below (due to clipping), which represents a singularity or the loss of mass, respectively) to orange (mass conservation: $\det(\mat{F}_1) = 1$) to white (expansion: $\det(\mat{F}_1) > 1$; white represents values of 2 or greater (due to clipping)). An illustration of this color map can be found in \figref{f:defgrid-jacobian-illustrated}. Notice that none of the maps for the determinant of the deformation gradient reported in this study values smaller than 0 (i.e., we do not report any degenerate deformation maps).

\section{Algorithm}
\label{s:algorithm}

Here, we provide more insight into our algorithm. We have added this information to the appendix, since we are mainly concerned with new regularization schemes in the present work. We refer to~\cite{Mang:2015a} for a detailed study of our globalized, inexact, preconditioned, reduced space (Gauss--)Newton--Krylov method for constrained diffeomorphic image registration; the study in \cite{Mang:2015a} includes a comparison to a preconditioned gradient descent scheme. Note that our solver is not finalized. We are currently working on improvements to reduce the time to solution.

\subsection{The Reduced Space Newton--Krylov Method}
\label{s:rsnkmethod}

We have seen in \secref{s:optimality-system} that the first order optimality conditions for~\eqref{e:opt-prob-v} are a system of space-time multicomponent nonlinear PDEs for the transported intensities $m$, the velocity field $\vect{v}$, and the mass source $w$. As we have seen in \secref{s:elimination} we can significantly simplify this system by exploiting variable elimination techniques to obtain~\eqref{e:first-order-opt-elim}; our algorithm will only operate on the reduced system. Note that this elimination not only reduces the computational complexity but also allows us to fulfill the hard constraint on $\idiv\vect{v}$ exactly. Given that we use a pseudospectral discretization in space, we can efficiently evaluate the resulting differential operators and their inverses.

We apply a Newton linearization to solve the first order optimality system~\eqref{e:first-order-opt}. This linearization results in a huge, severely ill-conditioned, dense, multi-component system for the incremental state, adjoint, and control variables. In \emph{full space methods} one directly solves this system. Our algorithm belongs to the class of \emph{reduced space methods} (see e.g.~\cite{Biros:2005a,Biros:2005b} for details); we do not solve for all unknowns of the KKT system simultaneously but eliminate the incremental state and adjoint variables from the system; we only iterate on the reduced space of the control variable $\vect{v}$. Advantages of reduced space methods include better spectral properties of the reduced space Hessian, a similar structure of the forward and adjoint operators, and a reduction of the order of the KKT system (the reduced space Hessian is nothing but the Schur complement for $\vect{\tilde{v}}$) to a size that is manageable~\cite{Biros:2005a}. Nonetheless, solving this system remains a significant challenge, given that the reduced space Hessian is still a large, ill-conditioned, dense, and compact operator.\footnote{A study of the spectral properties of the reduced space Hessian for compressible and incompressible diffeomorphisms can be found in~\cite{Mang:2015a}.} Further, as we will see below, we have to solve the state and adjoint equations at each iteration---a direct consequence of the block elimination in reduced space methods.

Next, we will discuss how the conceptual idea of a reduced space Newton--Krylov method relates to the optimality systems we have presented in \secref{s:optimality-conditions}. The reduced space KKT system is given by~\eqref{e:inc-control-eq-elim}. To solve this system, we have to evaluate to what we refer to as the reduced gradient $\vect{g}$ on the right-hand side of~\eqref{e:inc-control-eq-elim}. This involves the solution of the system~\eqref{e:first-order-opt-elim} in sequential order: Given some $\vect{v}$ we first solve~\eqref{e:state-eq-elim} forward in time. This gives us $m$ at $t=1$, which we need for the terminal condition in~\eqref{e:init-adj-eq-elim}.\footnote{We assign the costs for this forward solve to the evaluation of the objective and not the gradient.} Next, we solve~\eqref{e:adj-eq-elim} backward in time. Note that~\eqref{e:adj-eq-elim} represents a transport equation for the mismatch between $m_1$ and $m_R$. Thus, $\lambda$ will (ideally) tend to zero as we approach a solution (local minimizer) of~\eqref{e:opt-prob-v}. After we have solved~\eqref{e:state-eq-elim} and~\eqref{e:adj-eq-elim}, we can evaluate the control equation~\eqref{e:control-eq-elim} for some trial $\vect{v}$. Note that $m$ and $\lambda$ in~\eqref{e:control-eq-elim} are essentially functions of $\vect{v}$ through~\eqref{e:state-eq-elim} and~\eqref{e:adj-eq-elim}, respectively. Thus, although the PDE constraints are linear in $\vect{v}$, the inverse problem is not; it is non-linear in $\vect{v}$. This non-linearity as well as the conditioning of our problem are the main reasons why we prefer second order optimization methods.

Now, that we have found $\vect{g}$, we can solve the reduced space KKT system in~\eqref{e:inc-control-eq-elim} for the incremental control variable $\vect{\tilde{v}}$, i.e., compute the update (search direction) for $\vect{v}$; \eqref{e:inc-control-eq-elim} only provides the action of the reduced-space Hessian $\D{H}$ on $\vect{\tilde{v}}$ (Hessian matvec); this is all we need, given that we use a PCG method to solve~\eqref{e:inc-control-eq-elim} (i.e., our solver is matrix free). The incremental state and adjoint variables $\tilde{m}$ and $\tilde{\lambda}$ are---like we have seen for $m$ and $\lambda$ for the first order optimality conditions---functions of $\vect{\tilde{v}}$ through~\eqref{e:inc-state-eq-elim} and~\eqref{e:inc-adj-eq-elim}, respectively. Thus, each time we apply $\D{H}$ to $\vect{\tilde{v}}$, we have to solve~\eqref{e:inc-state-eq-elim} and~\eqref{e:inc-adj-eq-elim}. Also note that~\eqref{e:second-order-elim} not only depends on the incremental variables $\tilde{m}$, $\tilde{\lambda}$, and $\vect{\tilde{v}}$ but also on $m$, $\lambda$, and $\vect{v}$; the systems are strongly coupled.

To compute a minimizer to~\eqref{e:opt-prob-v} we have to repeat this entire process several times until convergence. We refer to the steps for updating $\vect{v}$ as \iname{outer iterations} and the steps for iteratively solving the reduced space KKT system as \iname{inner iterations}. We summarize these steps in compact form in \algref{a:outer-iteration} and \algref{a:inner-iteration} (for a PCG method), respectively.

The number of outer iterations depends on the rate of convergence of our scheme. A convergence study of our (Gauss--)Newton--Krylov scheme can be found in~\cite{Mang:2015a}. For a Newton--Krylov method we expect this convergence to be quadratic. However, given that we can not guarantee that the Hessian is positive definite far away from a (local) minimizer, we resort to a Gauss--Newton approximation. This is equivalent to dropping all expressions in which $\lambda$ appears in~\eqref{e:second-order-elim}; we expect the rate of convergence to drop from quadratic to super-linear (see~\cite{Mang:2015a} for more details). The number of inner iterations depends on the tolerance of the PCG method and on the preconditioner for the KKT system. For the tolerance we follow standard numerical optimization literature~\cite{Dembo:1983a, Eisenstat:1996a} and define it to be proportional to the relative $\ell^2$-norm of the reduced gradient (see~\cite{Mang:2015a} for details). The preconditioner is what we describe next.

\begin{algorithm}
\caption
{
Outer iteration of the designed inexact Newton--Krylov method.
}
\label{a:outer-iteration}
\algadjust
\begin{algorithmic}[1]
\STATE
{
$\vect{v}^h_0 \leftarrow 0$;
\quad
compute $m^h_0$, $\lambda^h_0$, $\F{J}^h(\vect{v}^h_0)$ and $\vect{g}^h_0$;
\quad
$k\leftarrow0$
}
\WHILE{true}
\STATE
    {
    stop $\leftarrow$ check for convergence
    }
    \STATE
    {
    \textbf{if} stop \textbf{break}
    }
    \STATE
    {
    $\vect{\tilde{v}}^h_k\leftarrow$
    solve~\eqref{e:reduced-space-kkt-system}
    given $m^h_k$, $\lambda^h_k$, $\vect{v}^h_k$, and $\vect{g}^h_k$
    \ccomment{Newton step (see~\algref{a:inner-iteration})}
    }
    \STATE
    {
    $\alpha_k \leftarrow$ perform line search on $\vect{\tilde{v}}^h_k$
    subject to the Armijo condition
    }
    \STATE
    {
    $
    \vect{v}^h_{k+1} \leftarrow \vect{v}^h_k
    + \alpha_k\vect{\tilde{v}}^h_k
    $
    }
    \STATE
    {
    $m^h_{k+1}(t=0) \leftarrow m_T^h$
    }
    \STATE
    {
    $m^h_{k+1} \leftarrow$ solve~\eqref{e:state-eq-elim} forward in time
    given $\vect{v}^h_{k+1}$
    \ccomment{forward solve}
    }
    \STATE
    {
    $\lambda^h_{k+1}(t=1) \leftarrow (m^h_R-m^h_{k+1}(t=1))$
    }
    \STATE
    {
    $\lambda^h_{k+1} \leftarrow$ solve~\eqref{e:adj-eq-elim}
    backward in time given $\vect{v}^h_{k+1}$ and $m^h_{k+1}$
    \ccomment{adjoint solve}
    }
    \STATE
    {
    compute $\F{J}^h(\vect{v}_{k+1}^h)$ and $\vect{g}^h_{k+1}$
    given $m^h_{k+1}$, $\lambda^h_{k+1}$ and $\vect{v}^h_{k+1}$
    }
    \STATE
    {
    $k \leftarrow k + 1$
    }
\ENDWHILE
\end{algorithmic}
\end{algorithm}

\begin{algorithm}
\caption
{Newton step. We illustrate the solution of the reduced KKT system~\eqref{e:reduced-space-kkt-system} using a PCG method at a given outer iteration $k\in\ns{N}$. The steps to compute the Hessian matrix vector product are given in lines~\ref{a:init-inc-adj-solve}--\ref{a:apply-hessian}.}
\label{a:inner-iteration}
\algadjust
\begin{algorithmic}[1]
\STATE
{
$\eta_k \leftarrow\min(0.5,(\|\vect{g}^h_k\|_2/\|\vect{g}^h_0\|_2)^{1/2})$
}
\STATE
{
$\vect{\tilde{v}}_0^h \leftarrow 0$,
\quad
$\vect{r}_0 \leftarrow - \vect{g}^h_k$,
\quad
$\vect{z}_0 \leftarrow (\F{A}^h)^{-1}\vect{r}_0$,
\quad
$\vect{s}_0 \leftarrow \vect{z}_0$,
\quad
$l\leftarrow 0$
}
\WHILE{$l < n$}
    \STATE
    {
    $\tilde{m}^h_l(t=0) \leftarrow 0$
    \label{a:init-inc-adj-solve}
    }
    \STATE
    {
    $\tilde{m}^h_l \leftarrow$ solve~\eqref{e:inc-state-eq-elim}
    forward in time
    given $m_k^h$, $\vect{v}^h_k$ and $\vect{\tilde{v}}^h_l$
    \ccomment{inc. forward solve}
    }
    \STATE{
    $\tilde{\lambda}^h_l(t=1) \leftarrow -\tilde{m}^h_l(t=1)$
    }
    \STATE
    {
    $\tilde{\lambda}^h_l \leftarrow$ solve~\eqref{e:inc-adj-eq-elim}
    backward in time
    given $\lambda_k^h$, $\vect{v}^h_k$ and $\vect{\tilde{v}}^h_l$
    \ccomment{inc. adjoint solve}
    }
    \STATE
    {
    $\vect{\tilde{s}}_l \leftarrow$ apply $\D{H}^h_l$
    to $\vect{s}_l$ as indicated in \eqref{e:inc-control-eq-elim}
    given $\lambda_k^h$, $\tilde{\lambda}_l^h$,
    $m^h_k$ and $\tilde{m}^h_l$
    \label{a:apply-hessian}
    }
    \STATE
    {
    $
    \kappa_l \leftarrow
    \langle\vect{r}_l,\vect{z}_l\rangle
    /\langle\vect{s}_l,\tilde{\vect{s}}_l\rangle
    $
    }
    \STATE
    {
    $
    \vect{\tilde{v}}^h_{l+1} \leftarrow \vect{\tilde{v}}^h_l
    + \kappa_l\vect{s}_l
    $
    }
    \STATE
    {
    $\vect{r}_{l+1}\leftarrow\vect{r}_l-\kappa_l\vect{\tilde{s}}_l$
    }
    \STATE
    {
    \textbf{if} $\|\vect{r}_{l+1}\|_2 < \eta_k$ \textbf{break}
    }
    \STATE
    {
    $\vect{z}_{l+1} \leftarrow (\D{A}^h)^{-1}\vect{r}_{l+1}$
    }
    \STATE
    {
    $
    \mu_l \leftarrow
    \langle\vect{z}_{l+1},\vect{r}_{l+1}\rangle
    /\langle\vect{z}_l,\vect{r}_l\rangle
    $
    }
    \STATE
    {
    $\vect{s}_{l+1} \leftarrow \vect{z}_{l+1} + \mu_l\vect{s}_l$
    }
    \STATE{$l \leftarrow l + 1$}
\ENDWHILE
\end{algorithmic}
\end{algorithm}

\subsection{Preconditioning the Reduced Space KKT System}

Preconditioning~\eqref{e:reduced-space-kkt-system} is essential to provide an efficient solver. We consider a left preconditioner $\mat{P}$ based on the second variation of the quadratic regularization models that act on $\vect{v}$; that is $\mat{P} \defeq \D{A}^h$.\footnote{By modifying the kernel of $\D{A}^h$, it is ensured that $\mat{P}$ is invertible.} This is a common choice in the PDE constrained optimization community. This preconditioner has---for our numerical scheme---essentially no construction and application cost; the inversion and application of $\D{A}^h$ amounts to a spectral diagonal scaling. The system PCG sees is a compact perturbation of the identity:\footnote{We slightly abuse our notation to indicate the dependence of $\D{L}[\vect{\tilde{b}}]$ on $\vect{\tilde{v}}$ trough \eqref{e:inc-state-eq-elim} and \eqref{e:inc-adj-eq-elim}, respectively; see \secref{s:rsnkmethod} for additional details.}
\begin{equation}
(\beta_v\D{A})^{-1} \left(\beta_v\D{A}[\vect{\tilde{v}}] + \D{L}[\vect{\tilde{b}}](\vect{\tilde{v}}) \right)
= \vect{\tilde{v}} + (\beta_v\D{A})^{-1}\D{L}[\vect{\tilde{b}}](\vect{\tilde{v}})
= (\id + (\beta_v\D{A})^{-1}\D{L}[\vect{\tilde{b}}])\vect{\tilde{v}};
\end{equation}

\noindent the inverse of the operator $\beta\D{A}$ acts like as smoother on $\D{L}[\vect{\tilde{b}}]$.

Since the second variation of~\eqref{e:shear-control-reg-v} can not directly be inverted in Fourier space (due to the complicated structure and the spatially varying viscosity; see~\eqref{e:2nd-var-reg-v-nls}) we use the inverse of the vectorial Laplacian operator as a preconditioner in case the regularization model in~\eqref{e:shear-control-reg-v} is considered.

\section{Performance Measures}
\label{s:measures-reg-performance}

We report different measures of registration performance. We summarize the definitions of these measures in \tabref{t:registration-performance-measures}. The inversion accuracy of our solver is controlled on the basis of a tolerance for the relative change of the reduced gradient. The quality of the inversion is assessed in terms of the relative change of the $L^2$-distance (residual) between the images to be registered. For some of our experiments we report values that measure the agreement between label maps of anatomical structures. We quantify regularity of the deformation map $\vect{y}$ on the basis of measures computed from the deformation gradient (see \secref{s:deformation-gradient}). In particular, we report values for the determinant of the deformation gradient. These indicate local volume change and deformation regularity. We also report values for the distance of the deformation gradient from identity.

\begin{table}
\caption{Summary of the performance measures. We report values for the relative change of the gradient (row 1) and the residual (row 2) to indicate inversion accuracy. Here, $\vect{g}^h_0$ is the initial gradient and $\vect{g}^h_{k^\star}$ is the gradient at the final iteration $k^\star$; $m_R^h$ is the reference and $m_T^h$ the template image. We also report measures derived from deformation gradient $\mat{F}_1^h$ (row 3 and row 4). Where applicable, we also report overlap measures (rows 5-8) for `'ground truth`' segmentations of anatomical structures. Here, $L_R$ denotes a label map for the reference image $m_R^h$ and $L_T$ the corresponding label map for the template image $m_T^h$; $\#$ denotes the cardinality of a set, $\backslash$ denotes the complement, and $\cup$ and $\cap$ are the union and intersection of two sets, respectively.}
\label{t:registration-performance-measures}
\centering
\tabadjust
\begin{tabular}{lcc}
\toprule
  Description
& Symbol
& Definition
\\\midrule
  relative change: reduced gradient
& $\|\vect{g}^{\star}\|_{\rm rel}$
& $\|\vect{g}^h_{k^\star}\|_2^2 / \|\vect{g}^h_0\|_2^2 $
\\
  relative change: residual
& $\|\vect{r}^{\star}\|_{\rm rel}$
& $\|m^h_R - m^h_1\|_2^2 / \|m^h_R - m^h_T\|_2^2 $
\\
  distance of deformation gradient from identity
& $D$
& $\|\mat{F}_1^h - \mat{I}\|_F$
\\
  determinant of the deformation gradient
& $J$
& $\det(\mat{F}_1^h)$
\\
  overlap: Jaccard similarity coefficient
& JSC
& $\# (L_R \cap  L_T)/ \#(L_R \cup  L_T)$
\\
  overlap: Dice similarity coefficient
& DSC
& $ 2\#(L_R \cap  L_T)/(\# L_R + \# L_T)$
\\
  overlap: false positive error
& FPE
& $\#(L_T \backslash  L_R)/\# L_T$
\\
  overlap: false negative error
& FNE
& $\#(L_R \backslash  L_T)/\# L_R$
\\\bottomrule
\end{tabular}
\end{table}

\end{appendix}

\end{document}